\documentclass[a4paper,oneside,english,10pt]{article}
\usepackage[T1]{fontenc}
\usepackage[latin9]{inputenc}
\usepackage{color}
\usepackage{amsthm}
\usepackage{amsmath}
\usepackage{amstext}
\usepackage{graphicx}
\usepackage{float}
\usepackage{amssymb}
\usepackage{esint}
\usepackage{a4wide}
\usepackage{babel}
\usepackage{bbm}
\usepackage{mathrsfs}
\usepackage{upgreek}
\usepackage{setspace} 
\usepackage{multirow}
\usepackage{placeins}
\usepackage{array}
\usepackage{hyperref} 
\usepackage{algorithmic}
\usepackage{algorithm}

\usepackage{setspace}

\newtheorem{theorem}{Theorem}[section]

\newtheorem{definition}[theorem]{Definition}
\newtheorem{remark}[theorem]{Remark}

  \newcommand{\vnabla}{\boldsymbol{\nabla  }}

\newcommand{\Rot}{\operatorname{\boldsymbol{curl}}}
\newcommand{\Div}{\operatorname{\mathrm{div}}}

\newcommand  {\vphi}{\boldsymbol{\varphi}}

\newcommand{\uinc} {\uu^{\mathrm{inc}}}

\newcommand  {\R}{{\mathbb R}}
 
\newcommand  {\nn}{\boldsymbol{n}}
\newcommand  {\uu}{\boldsymbol{ u}}
  \newcommand {\Id}{\mathbf{I}} 
   \newcommand{\transposee}[1]{{\vphantom{#1}}{#1}^{\sf T}}  
\newcommand{\vpsi}{{\boldsymbol{\psi}}}
\newcommand  {\HH}{\boldsymbol{ H}}
\newcommand  {\TT}{\bf{T}}
\newcommand  {\LL}{\boldsymbol{L}}
 \newcommand  {\x}{\boldsymbol{x}}
  \newcommand  {\y}{\boldsymbol{y}}
    \renewcommand  {\tt}{\boldsymbol{ t}}
    
        \usepackage{authblk}
    
\begin{document}

\title{Theory and implementation of $\mathcal{H}$-matrix based iterative and direct solvers for  
 {Helmholtz and elastodynamic}
oscillatory kernels}

\author[1]{St\'ephanie Chaillat\thanks{stephanie.chaillat@ensta-paristech.fr}}

\author[1]{Luca Desiderio\thanks{luca.desiderio@yahoo.it}}

\author[1]{Patrick Ciarlet\thanks{patrick.ciarlet@ensta-paristech.fr}}

\affil[1]{Laboratoire POEMS UMR CNRS-INRIA-ENSTA, Université Paris-Saclay\\ ENSTA-UMA,  828 Bd des Mar\'echaux, 91762 Palaiseau Cedex, FRANCE}%

\maketitle

\begin{abstract} 

In this work, we study  the accuracy and efficiency of hierarchical matrix ($\mathcal{H}$-matrix) based fast methods for solving dense linear systems arising from the discretization of the 3D elastodynamic Green's tensors. It is well known in the literature that standard  $\mathcal{H}$-matrix based methods,  although very efficient tools for asymptotically smooth kernels, are  not optimal for oscillatory kernels. $\mathcal{H}^2$-matrix and  directional approaches have been proposed to overcome this problem. However the implementation of such methods is much more involved than the standard $\mathcal{H}$-matrix representation. 
The central questions we address  are twofold. (i) What is the frequency-range in which the  $\mathcal{H}$-matrix format is an efficient representation for 3D elastodynamic problems? (ii) What can be expected of such an approach to model problems in mechanical engineering?
 We show   that even though the method is not optimal (in the sense that more involved representations can lead to faster algorithms) an efficient solver can be easily developed. The capabilities of the method are illustrated on numerical examples using the Boundary Element Method.
\end{abstract}

Time-harmonic elastic waves; $\mathcal{H}$-matrices, Low-rank approximations; Estimators; Algorithmic complexity; fast BEMs.

\section{Introduction}

The 3D linear isotropic elastodynamic equation for the displacement field $\uu$ (also called Navier equation) is given by
\begin{equation}
\label{NE}
\Div\upsigma(\uu)+ \rho\omega^2\uu=0
\end{equation}
where $\omega>0$ is the circular  frequency. It is supplemented with appropriate boundary conditions which contain the data. 
The stress and strain tensors are respectively given by 
$\upvarsigma(\uu)=\lambda(\Div\uu)\Id_{3}+2\mu\upvarepsilon(\uu)$ and $\upvarepsilon(\uu)=\dfrac{1}{2}\big([\vnabla\uu]+\transposee{[\vnabla\uu]}\big)$, where $\Id_{3}$ is the 3-by-3 identity matrix  and $[\vnabla\uu]$ is the 3-by-3 matrix whose $\beta$-th column is the gradient of the $\beta$-th component of $\uu$ for $1\le \beta \le 3$,  
  $\mu$ and $\lambda$ are the Lamé parameters  and  $\rho$ the density.  
Denoting $\kappa_{p}^2=\rho\omega^2(\lambda+2\mu)^{-1}$ and $\kappa_{s}^2=\rho\omega^2\mu^{-1}$  the so-called P and S wavenumbers,
 the Green's tensor of the Navier equation is  a 3-by-3 matrix-valued function expressed by
\begin{equation}
{\bf U}_{\omega}(\x,\y)=\dfrac{1}{\rho\omega^2}\left(\Rot\Rot_{\x} \left[\dfrac{e^{i\kappa_{s}|\x-\y|}}{4\pi|\x-\y|}\,{ \Id}_{3}\right]-\vnabla_{\x}\boldsymbol{ \Div}_{\x}\,\left[\dfrac{e^{i\kappa_{p}|\x-\y|}}{4\pi|\x-\y|}\Id_{3}\right]\right)
\label{elasto_U}
\end{equation}
where the index $\x$ means that differentiation is carried out with respect to $\x$ and  $\boldsymbol{ \Div_x} \mathbb{A}$   corresponds to the application of the divergence along each row of $\mathbb{A}$. 
One may use  this tensor to  represent the solution of~(\ref{NE}). Alternately,  one may use the 
  tensor $\TT_{\omega}(\x,\y)$, which is obtained by applying the traction operator 
 \begin{equation}
 {\boldsymbol T}=2\mu\dfrac{\partial}{\partial\nn}+\lambda\nn\Div+\mu\,\nn\times\Rot
 \label{elasto_T}
\end{equation}
 to each column of $\bf{U}_{\omega}(\x,\y)$: $\TT_{\omega}(\x,\y)=[ {\boldsymbol T}_{\y} {\bf U}_{\omega}(\x,\y)]$ where the index $\y$ means that differentiation is carried out with respect to $\y$. 
 
We consider   the fast solution of dense linear systems   of the form
\begin{equation}
\mathbb{ A} {\bf p}={\bf b}, \quad {\mathbb A}\in \mathbb{C}^{3N_c \times 3N_c}
\label{syst_a}
\end{equation}
where ${\mathbb A}$ is the matrix corresponding to the discretization of the 3-by-3 Green's tensors ${\bf U}_{\omega}(\x_i,\y_j)$ or ${\bf T}_{\omega}(\x_i,\y_j)$ for two clouds of $N_c$ points $ {(\x_i)_{ 1 \le i \le N_c}}$ and $ {(\y_j)_{1 \le j \le N_c}}$. Here ${\bf p}$ is the unknown vector approximating the solution at $ {(\x_i)_{ 1 \le i \le N_c}}$ and ${\bf b}$ is a given right hand side that depends on the data. Such dense systems are encountered for example in the context of the Boundary Element Method~\cite{bonnet1999boundary,sauter2010boundary}.

  If no compression or acceleration technique is used, the storage of such a system is of the order $O(N_c^2)$, the iterative solution (e.g. with GMRES) is $O(N_{\operatorname{iter}}N_c^2)$ where $N_{\operatorname{iter}}$ is the number of iterations, while the direct solution (e.g. via LU factorizations)
 is  $O(N_c^3)$. In the last decades, different  approaches have been proposed to speed up the solution of  dense systems. 
  The most known  method is probably the fast multipole method  (FMM) proposed by Greengard and Rokhlin~\cite{greengard1987fast} which enables a fast evaluation of the matrix-vector products. We recall that the matrix-vector product is the crucial tool in the context of an iterative solution. Initially   developed for N-body simulations, the FMM has then been extended to oscillatory kernels~\cite{greengard1998accelerating,darve2000fast}. The method is now widely used in many application fields and has shown its capabilities in the context of mechanical engineering problems solved with the BEM~\cite{chaillat2013recent,takahashi2012wideband}. 
  
An alternative approach designed for dense systems is based on the concept of hierarchical matrices ($\mathcal{H}$-matrices)~\cite{bebendorf2008hierarchical}. The principle of $\mathcal{H}$-matrices is to partition the initial dense linear system, and then approximate it into a data-sparse one, by finding sub-blocks in the matrix that can be accurately estimated by low-rank matrices. In other terms, one further approximates the matrix $\mathbb{A}$ from~(\ref{syst_a}). The efficiency of hierarchical matrices relies on the possibility to approximate, under certain conditions, the underlying kernel function by low-rank matrices. 
The approach has been shown to be  very efficient for asymptotically smooth kernels (e.g. Laplace kernel).
On the other hand, oscillatory kernels such as the Helmholtz or elastodynamic kernels, are  not asymptotically smooth. In these cases, the method is not optimal~\cite{banjai2008hierarchical}.
To avoid the increase of the rank for high-frequency problems, 
directional  $\mathcal{H}^2$-methods have been proposed~\cite{borm2015approximation,borm2015directional}. $\mathcal{H}^2$-matrices are a specialization of hierarchical matrices. It is a multigrid-like version of $\mathcal{H}$-matrices  that enables more compression, by factorizing some basis functions of the approximate separable expansion~\cite{borm2006matrix}.  

Since the implementation of $\mathcal{H}^2$ or directional methods  is much more involved than the one of the standard $\mathcal{H}$-matrix, it is important to determine  
 the frequency-range within which the  $\mathcal{H}$-matrices are efficient for elastodynamic problems and  what    can be expected of such an approach   to solve problems encountered in mechanical engineering. Previous works on $\mathcal{H}$-matrices for oscillatory kernels have mainly   be  devoted to the direct application   to derive fast iterative solvers for 3D acoustics~\cite{brunner2010comparison, stolper2004computing}, a direct solver for 3D electromagnetism~\cite{lize} or iterative solvers for 3D elastodynamics~\cite{milazzo2012hierarchical,MESSNER2010}.
 There is no discussion in these references on the capabilities and limits of the method for oscillatory kernels.
 We show in this work that even though the method is not optimal (in the sense that more efficient approaches can be proposed at the cost  of a much more complex implementation effort), an efficient solver is easily developed. The capabilities of the method are illustrated on numerical examples using the Boundary Element Method.

 The article is organized as follows. After reviewing general facts about  $\mathcal{H}$-matrices in Section~2, 
  $\mathcal{H}$-matrix based solvers and an   estimator to certify the results of the direct solver are given in Section~3.
 Section~4 gives an overview of the  Boundary Element Method for 3D elastodynamics.
In Section~5, theoretical estimates for the application of $\mathcal{H}$-matrices  to  elastodynamics are derived. 
Section~6 presents the implementation issues for elastodynamics. Finally in Sections~7 and~8, we present numerical tests with varying number of points $N_c$ for both representations of the solution ($\bf{U}_{\omega}$ or $\bf{T}_{\omega}$). 
In Section~7, numerical tests for the low  frequency regime i.e. for a fixed frequency  are reported and discussed. In Section~8, similar results are given for  the high frequency regime i.e. for   a fixed density of points per S-wavelength.

\section{General description of the $\mathcal{H}$-LU factorization}

 \subsection{$\mathcal{H}$- matrix representation \label{section_tree}}
 
 Hierarchical matrices or $\mathcal{H}$-matrices have been introduced by Hackbusch~\cite{hackbusch1999sparse} to compute a data-sparse representation 
 of some special dense matrices (e.g. matrices resulting from the discretization of non-local operators). The principle of $\mathcal{H}$-matrices is (i) to partition  the matrix into blocks and (ii) to perform  low-rank approximations of  the blocks of the matrix which are known \emph{a priori} (by using an admissibility condition) to be accurately approximated by low-rank decompositions. With these two ingredients it is possible to define fast iterative and direct solvers for matrices having a hierarchical representation. 
Using low-rank representations, the memory requirements and costs of a matrix-vector product are reduced. In addition, using   $\mathcal{H}$-matrix arithmetic it is possible to derive fast direct solvers.

\begin{figure}[!htb]
 \begin{center}
   \begin{tabular}{cc}
 \includegraphics[width=0.45 \textwidth]{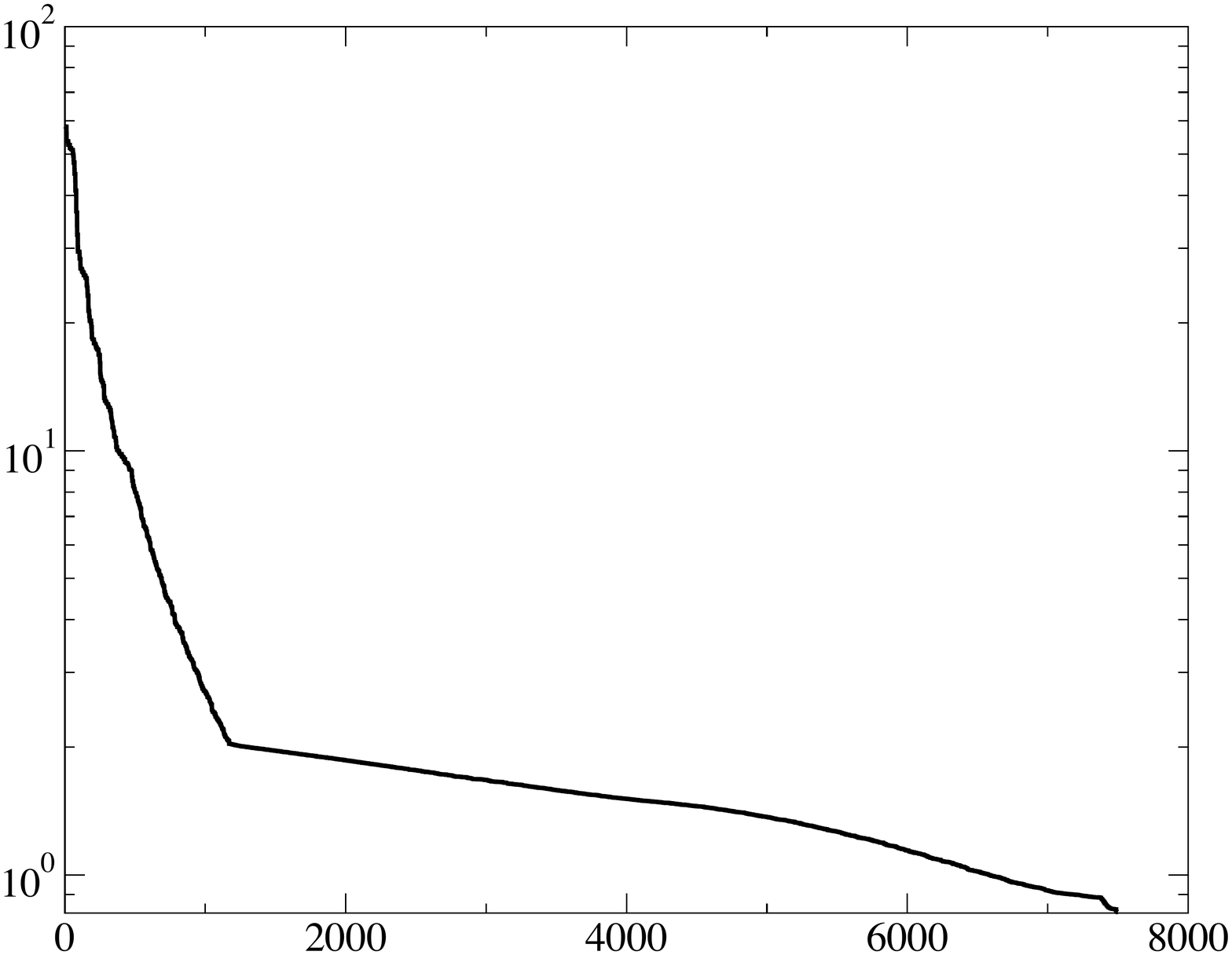} &  \includegraphics[width=0.45 \textwidth]{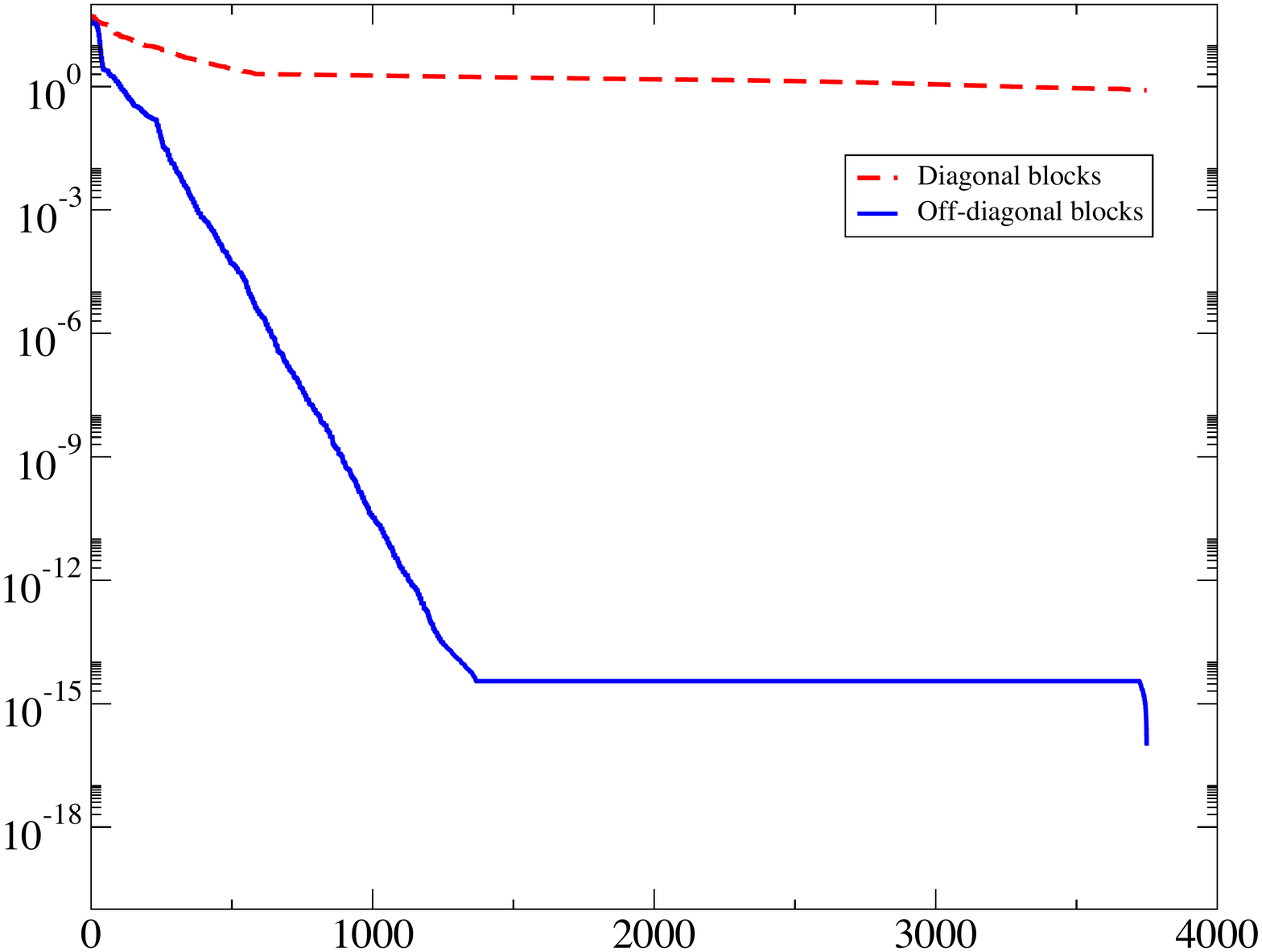}\\
 (a) & (b)
 \end{tabular}
 \end{center}
 \caption{(a) Decrease of the singular values of the complete matrix ${\mathbb G}_e$ corresponding to the discretization of the elastic Green's tensor for a cloud of $2 \ 500$ points. (b) Decrease of the singular values of the diagonal and off-diagonal blocks of the same matrix decomposed into a 2-by-2 block matrix. }
 \label{svd_all}
 \end{figure}

\begin{figure}[!htb]
 \begin{center}
\includegraphics[scale=0.5]{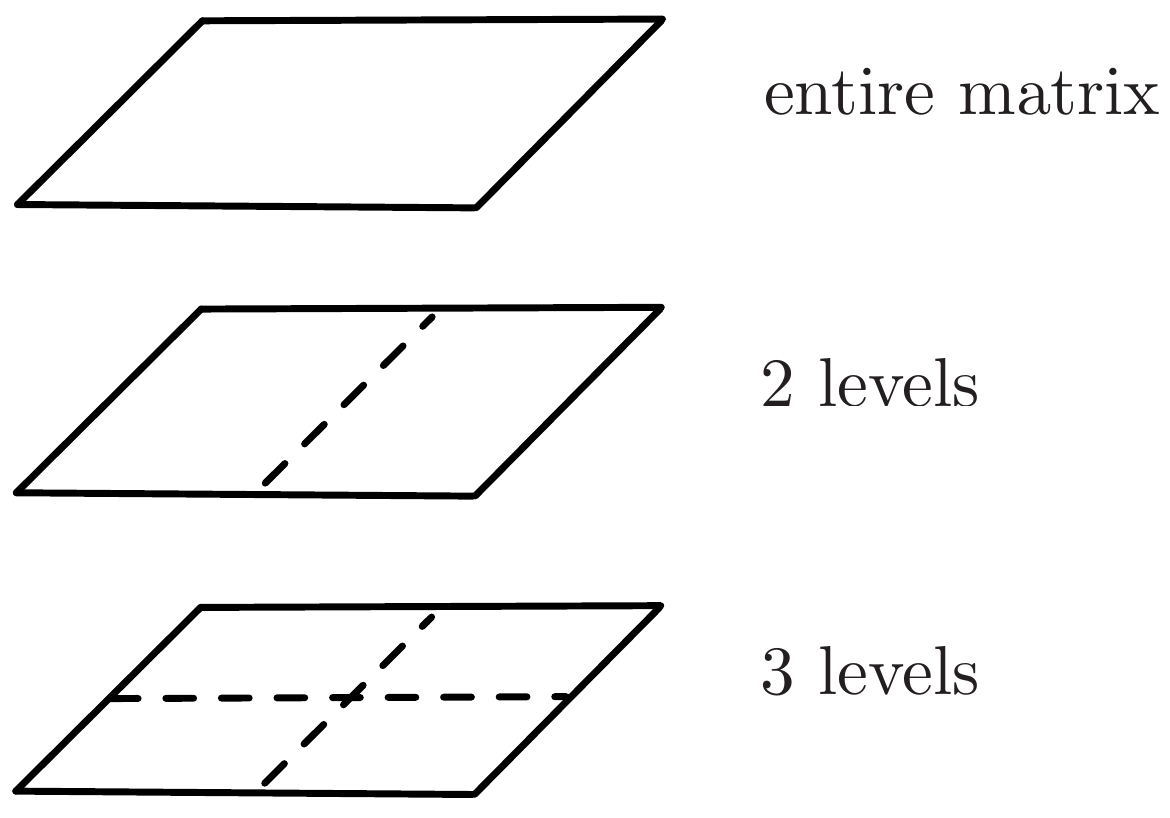}   
 \end{center}
 \caption{Illustrative example:  partition of the degrees of freedom. }
  \label{geo_illu}
 \end{figure}

 \paragraph{Illustrative example} To illustrate the construction of an $\mathcal{H}$-matrix, we consider the matrix ${\mathbb G}_e$ resulting from the discretization of the  3D elastic Green's tensor ${\bf U}_{\omega}$    for a cloud of points located on the plate $(x_1,x_2,x_3)\in [-a,a]\times[-a,a]\times \{0\}$. 
 The plate is discretized uniformly with a fixed density of $10$ points per S-wavelength $\lambda_s=2\pi/\kappa_s$.
 We fix the Poisson's ratio $\nu$ to $\nu=\frac{\lambda}{2(\lambda+ \mu)}=1/3$ and a non-dimensional frequency $\eta_S=\kappa_s a=5 \pi$ (e.g. $\rho\mu^{-1}=1$,  $\omega=5\pi$, $\kappa_s=5 \pi$, $\kappa_p=\kappa_s/2$  and $a=1$). 
  As a result,  the discretization consists of  $N_d=50$ points in each direction leading to $N_c=2500$ points and a matrix of size $7500 \times 7500$.  %
  We recall that the numerical rank of a matrix ${\mathbb A}$ is
 \begin{equation}
 r(\varepsilon):=\mbox{min} \{ r \quad | \quad   ||{\mathbb A} -{\mathbb A}_r|| \le \varepsilon ||{\mathbb A}|| \}
 \label{num_rank}
 \end{equation}
 where ${\mathbb A}_r$ defines the  singular value decomposition (SVD) of ${\mathbb A}$ keeping only the  $r$ largest singular values and $\varepsilon>0$ is a given parameter.
In Fig.~\ref{svd_all}a, we report the decrease of the singular values of the matrix ${\mathbb G}_e$. As expected, the decay of the singular values is very slow such that the matrix cannot be approximated by a low-rank decomposition.
 Now, we partition the plate into two parts of equal area (see Fig.~\ref{geo_illu}) and 
  subdivide the matrix accordingly into 4 subblocks. Figure~\ref{svd_all}b gives the decrease of the singular values both for diagonal and off-diagonal blocks. It illustrates the possibility to accurately represent off-diagonal blocks by low-rank matrices while diagonal blocks cannot have low-rank representations.
 
\begin{figure}[!htb]
 \begin{center}
\begin{tabular}{ccc}
  \includegraphics[scale=0.5]{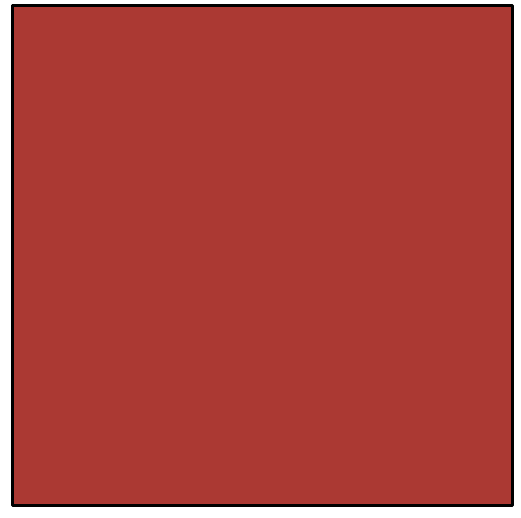} &\includegraphics[scale=0.5]{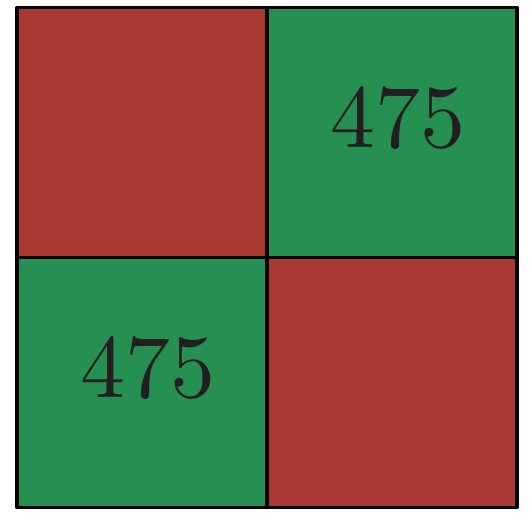} & \includegraphics[scale=0.5]{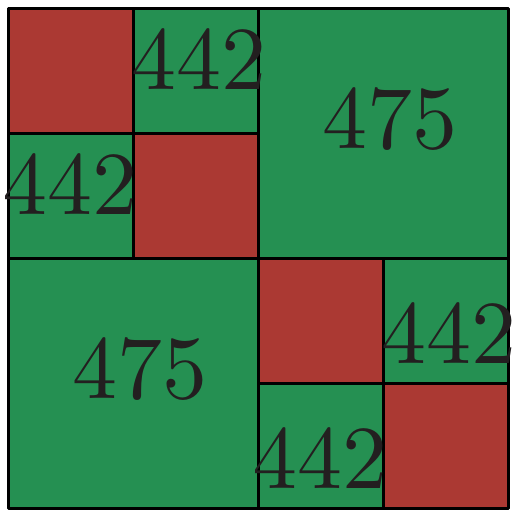} \\
  {Entire matrix} & {2 levels} &    {3 levels}\\ 
Size block=$7500$ & Size block=$3750$ & Size block=$1875$\\[-2ex]
\end{tabular}
\end{center}
 \caption{ Computed numerical ranks of each block of the matrix ${\mathbb G}_e$ to achieve an accuracy of $10^{-4}$, i.e. if singular values smaller than $10^{-4}$ are neglected in the singular value decomposition.}
  \label{ranks}
 \end{figure}

\begin{figure}[!htb]
 \begin{center}
 \begin{tabular}{cc}
\includegraphics[scale=0.65]{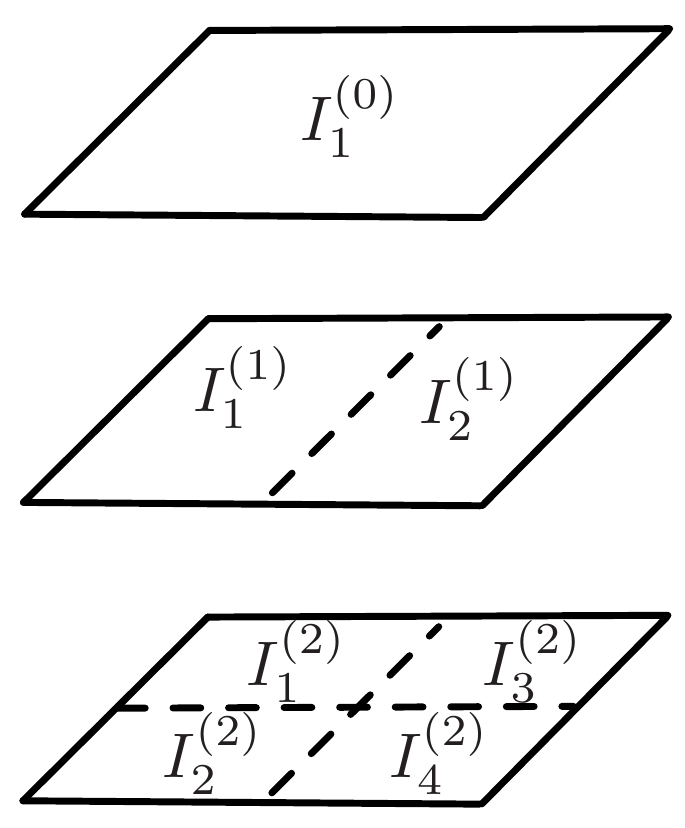} & \includegraphics[scale=0.35]{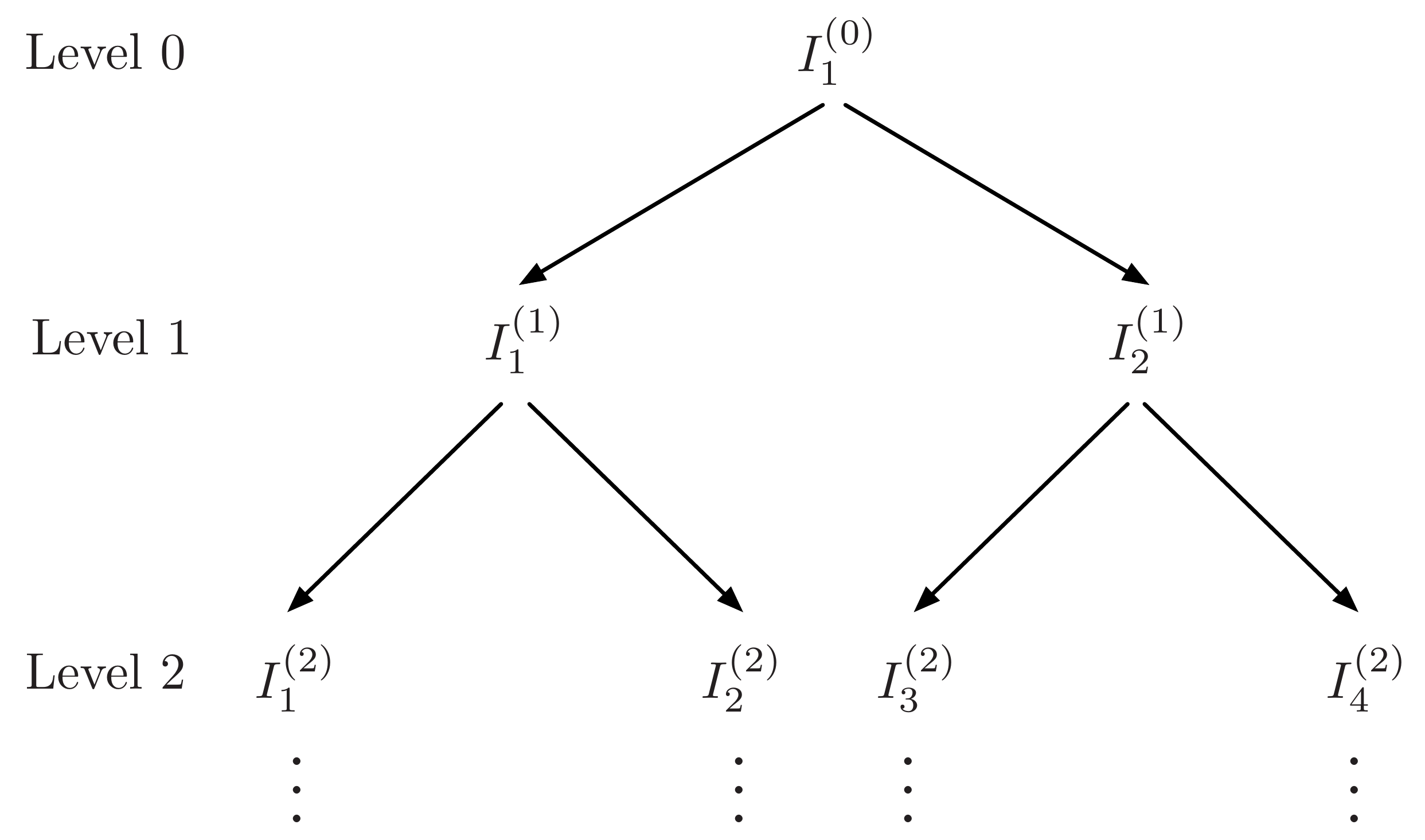} \\
 (a) Partition of the physical domain & (b) Binary cluster tree $\mathcal{T}_I$
 \end{tabular}
 \end{center}
  \caption{Illustration of the clustering of the degrees of freedom: (a) partition of the degrees of freedom in the domain and (b) corresponding binary tree.  }
 \label{binary_tree}
 \end{figure}

  If we keep subdividing the full rank blocks in a similar manner, we observe   in Fig.~\ref{ranks} that diagonal blocks are always full rank  numerically, i.e. with respect to~(\ref{num_rank})  where we have chosen $\varepsilon=10^{-4}$, while  off-diagonal blocks become accurately approximated by a low-rank decomposition after some iterations of the subdivision process. 
  This academic example illustrates the concept of data-sparse matrices used to derive fast algorithms. Using the hierarchical structure in addition to low-rank approximations, significant savings are obtained both in terms of computational times and memory requirements. 
 \paragraph{Clustering of the unknowns} 
 The key ingredient of hierarchical matrices is the recursive block subdivision.
 On the illustrative example, the subdivision of the matrix is conducted by a recursive subdivision of the geometry, i.e. a recursive subdivision of the plan   into partitions of equal areas. 
 The first step prior to the partition of the matrix is thus   a partitioning  based on the geometry of the set of   row and column indices of the matrix ${\mathbb A}$. 
 The purpose   is to permute the indices in the matrix   to reflect the physical distance and thus interaction between degrees of freedom. Consecutive indices should correspond to DOFs that interact at close range. For the sake of clarity,  in this work ${\mathbb A}$ is defined by the same set of indices  $I=\{1, \ldots, n\}$ for rows and columns. A binary tree $\mathcal{T}_I$ is used to drive the clustering. Each node of the tree defines a subset of indices $\sigma \subset I$
and each subset corresponds to a  
part in the partition of the domain, see Figure~\ref{binary_tree} for the case of the illustrative example.
There exist different  approaches to perform the subdivision~\cite{hackbusch2015hierarchical}.
We consider the simplest possible one~: based on a geometric argument. For each node in the tree, we determine the  box enclosing all the points in the cloud and subdivide it into 2 boxes, along the largest dimension.  The subdivision is stopped when a minimum number of DOFs per box is reached ($N_{\operatorname{leaf}}=100$ in the following). For uniform meshes, this strategy defines a balanced binary tree~\cite{bebendorf2008hierarchical} such that
 the number of levels in the tree $\mathcal{T}_I$ is given by $L(I)= \lceil  \log_2(\frac{|I|}{N_{\operatorname{leaf}}}) \rceil \le \log_2 |I| -\log_2  N_{\operatorname{leaf}} +1$.
 Note that for this step  only geometrical information is needed.

\begin{figure}[!htb]
 \begin{center}
 \begin{tabular}{cc}
 \includegraphics[width=7cm]{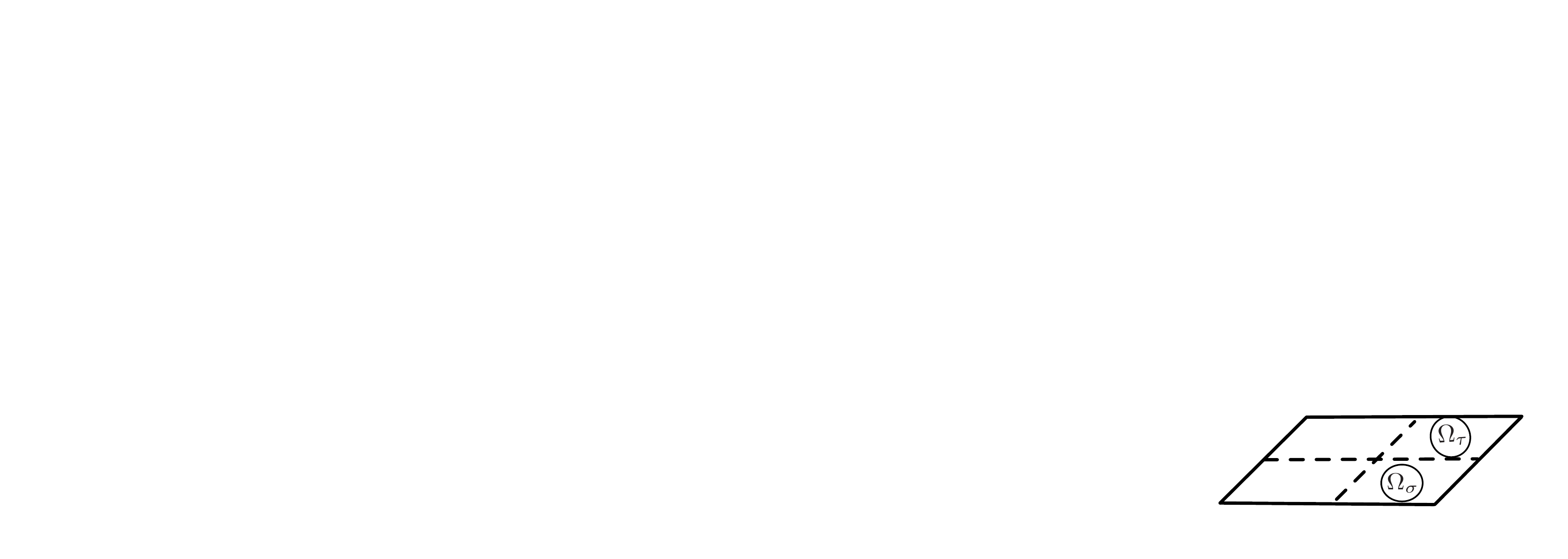} &  \includegraphics[width=4cm]{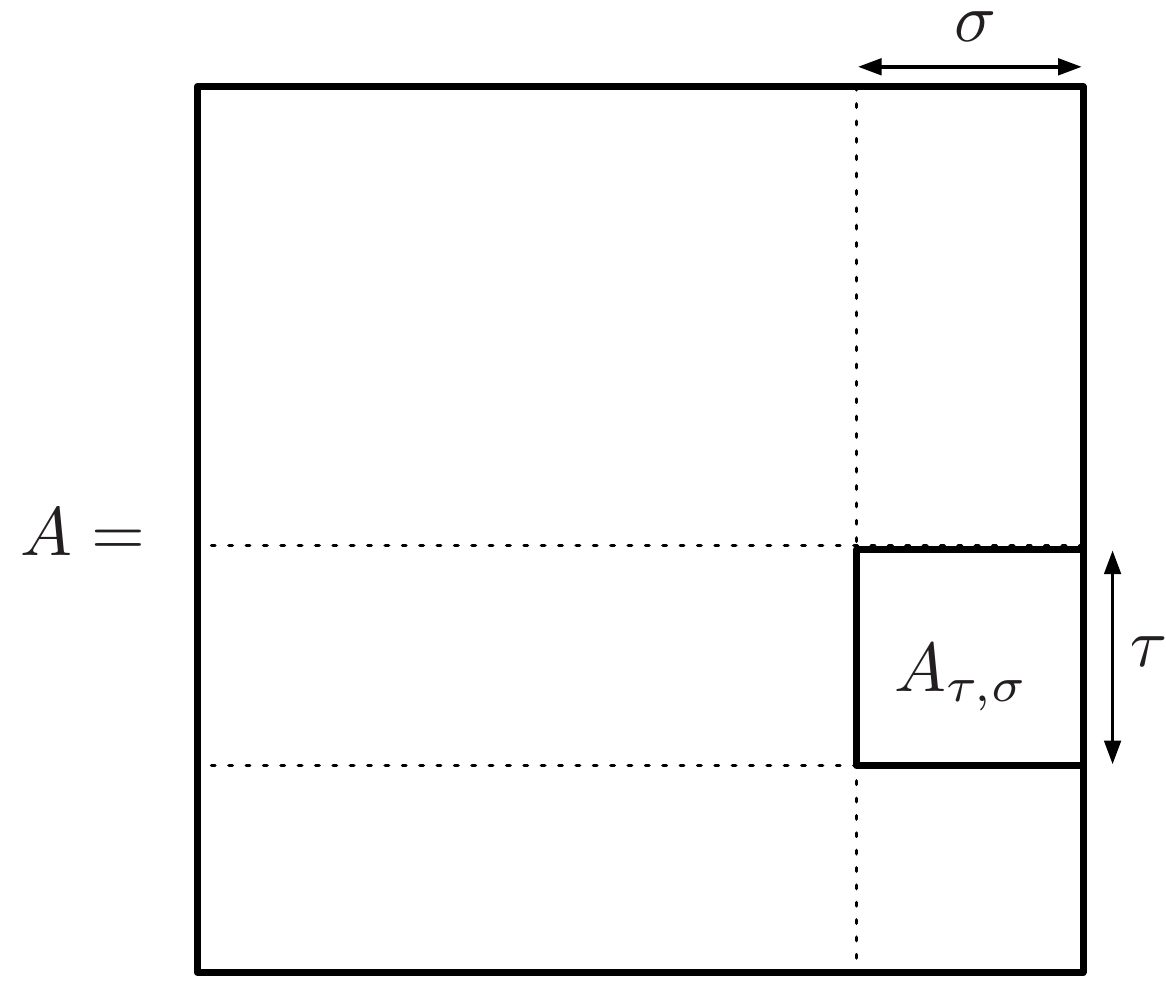} \\
 (a) & (b)
 \end{tabular}
 \end{center}
  \caption{Illustration of the construction of the block cluster tree: (a) Clustering of the unknowns on the geometry and (b) corresponding block clustering in the matrix. }
 \label{def_tree}
 \end{figure}

\begin{figure}[!htb]
 \begin{center}
 \begin{tabular}{cc}
  \includegraphics[width=4cm]{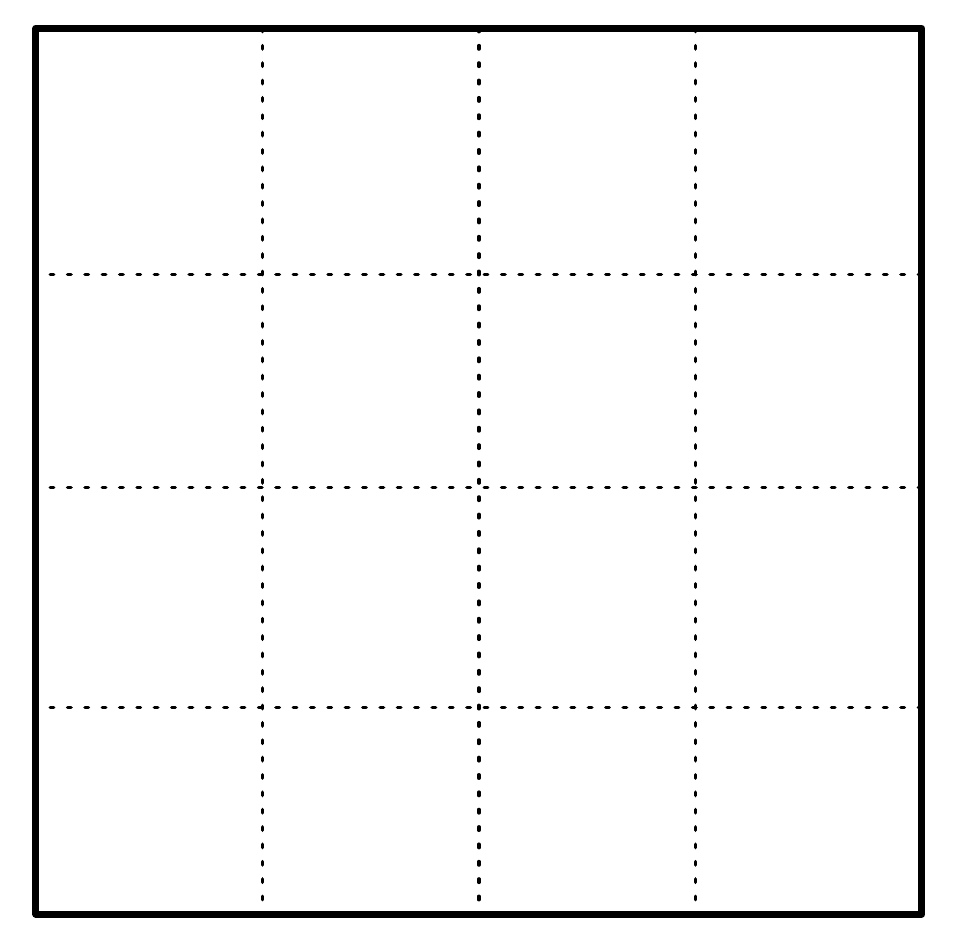} &  \includegraphics[width=4cm]{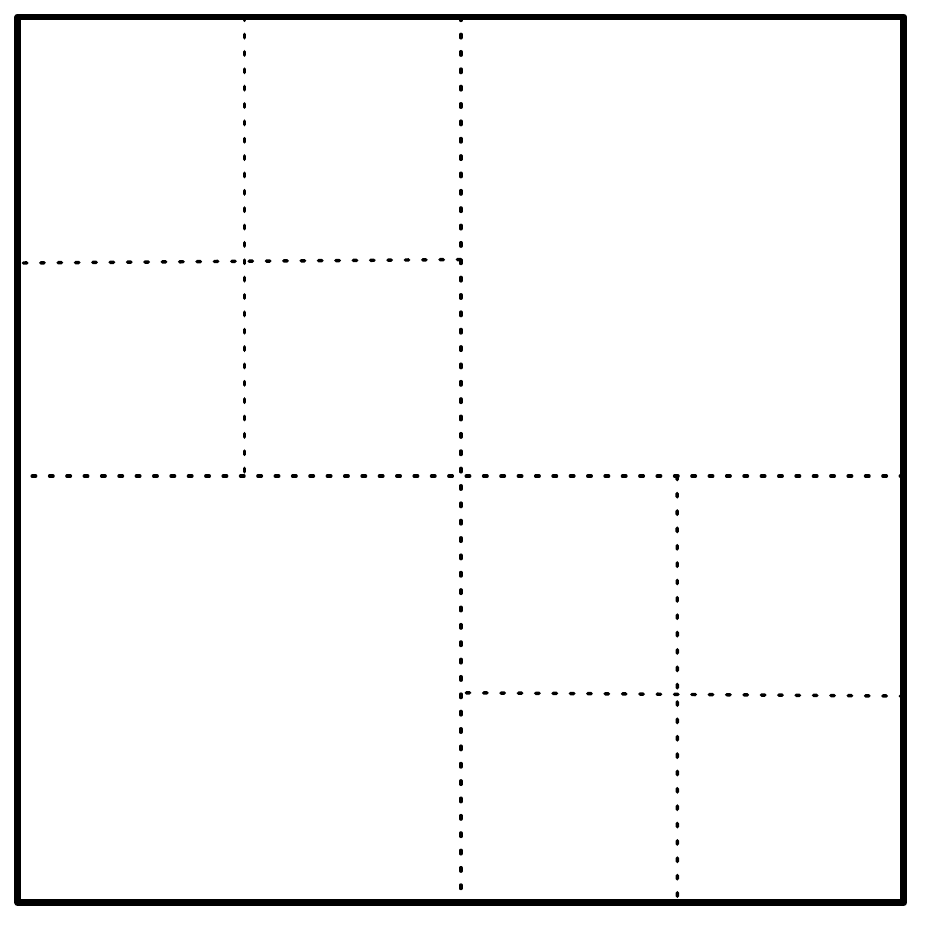}  \\
 (a) & (b)
 \end{tabular}
 \end{center}
  \caption{ (a) Block cluster representation $\mathcal{T}_{I \times I}$ for the illustrative example  (full structure); (b)  Hierarchical partition $\mathcal{P} \subset \mathcal{T}_{I \times I}$  of the same matrix based on the admissibility condition (sparse structure). }
 \label{def_tree_uni}
 \end{figure}

 \paragraph{Subdivision of the matrix} After the clustering of the unknowns is performed, a
   block cluster  representation  $\mathcal{T}_{I\times I}$   of the matrix ${\mathbb A}$ is defined  by going through the cluster tree $\mathcal{T}_{I}$. 
 Each node of  $\mathcal{T}_{I\times I}$ contains a pair $(\sigma,\tau)$ of indices of $\mathcal{T}_{I}$ and defines a block of ${\mathbb A}$  (see Figure~\ref{def_tree}). 
 This uniform partition defines a block structure of the matrix with a full pattern of $4^{L(I)-1}$ blocks, in particular every node of the tree at the leaf level is connected with all the other nodes at the leaf level (Figure~\ref{def_tree_uni}a).
  Figure~\ref{svd_all}b shows that this partition is not optimal. As a matter of fact, some  parts of the matrix ${\mathbb A}$ can accurately be  approximated by a low-rank matrix at a higher level (i.e. for larger clusters). Such blocks are said to be \emph{admissible}.  A hierarchical representation $\mathcal{P}\subset \mathcal{T}_{I \times I}$ that uses the cluster tree $\mathcal{T}_{I}$ and the existence of \emph{admissible} blocks is more appropriate (Figure~\ref{def_tree_uni}b). 
 Starting from the initial matrix, each block is recursively subdivided until it is either \emph{admissible} or the leaf level is achieved.
In the illustrative example, only diagonal blocks are subdivided due to the very simple 2D geometry (see again Fig.~\ref{ranks}, red-colored blocks). For more complex 3D geometries, an admissibility condition based on the geometry and the interaction distance between points is used to determine \emph{a priori} the \emph{admissible} blocks.
For more details on the construction of the block cluster tree, we refer the interested reader to~\cite{borm2003introduction}.
 The partition $\mathcal{P}$ is subdivided into two subsets $\mathcal{P}^{\operatorname{ad}}$  and $\mathcal{P}^{\operatorname{non-ad}}$ reflecting the possibility for a block $\tau \times \sigma$ to be either \emph{admissible}, i.e.  $\tau \times \sigma \in \mathcal{P}^{\operatorname{ad}}$; or \emph{non-admissible}, i.e.  $\tau \times \sigma \in \mathcal{P}^{\operatorname{non-ad}}$. It is clear that  $\mathcal{P}=\mathcal{P}^{\operatorname{ad}}\cup \mathcal{P}^{\operatorname{non-ad}}$.  To sum up, the blocks of the partition can be of 3 types: at the leaf level a block can be either an \emph{admissible} block  or a  \emph{non-admissible} block, at a non-leaf level  a block can be either   an \emph{admissible} block or an $\mathcal{H}$-matrix (i.e a block that will be subsequently hierarchically subdivided).
 
  At this point, we need to define the sparsity pattern introduced by L. Grasedyck~\cite{grasedyck2001theorie}. For sparse matrices, the sparsity pattern gives the maximum number of non-zero entries per row. Similarly for a hierarchical matrix defined by the partition $\mathcal{P} \subset \mathcal{T}_{I\times I}$, the sparsity pattern of a row cluster $\tau \in \mathcal{T}_I$, resp. a column cluster $\sigma \in \mathcal{T}_I$, is 
 \[
 C_{sp}(\tau)=|\{ \sigma \in \mathcal{T}_I: \tau \times \sigma \in \mathcal{P}\}|, \mbox{ resp. }  C_{sp}(\sigma)=|\{ \tau \in \mathcal{T}_I: \tau \times \sigma \in \mathcal{P}\}|. 
 \]
 It is convenient to define the overall sparsity pattern (for row and column clusters):
 \[
\displaystyle{  C_{sp}=\mbox{max} \{ \mbox{max}_{\tau \in \mathcal{T}_I} C_{sp}(\tau), \ \mbox{max}_{\sigma \in \mathcal{T}_I} C_{sp}(\sigma) \}.}
 \]
 
  \paragraph{Special case of asymptotically smooth kernels}
 
  $\mathcal{H}$-matrix representations  have been derived for some specific problems and will not result in  efficient algorithms    for all equations or matrices. The crucial point is to know 
  \emph{a priori} (i) if savings will be obtained when trying to approximate \emph{admissible} blocks with a sum of separated variable functions and (ii) which blocks are \emph{admissible} since the explicit computation of the rank of all the blocks would be too expensive.   In the case of asymptotically smooth kernels
  $G(\x,\y)$, it is proved that  under some \emph{a priori} condition on the distance between $\x$ and $\y$, the kernel is a degenerate function (see Section~\ref{theory}).
  After discretization, this property is reformulated as the efficient approximation of blocks of the matrix by low-rank  matrices. 
  The Laplace Green's function is an example of asymptotically smooth kernel for which $\mathcal{H}$-matrix representations have been shown to be very efficient.
Since the 
  3D elastodynamics Green's tensor, similarly to the Helmholtz Green's function, is   a linear combination of  derivatives of the 3D Laplace Green's function  with coefficients depending on the circular frequency $\omega$, this work is concerned
  with the determination of the frequency range for which  hierarchical  representations can be successful for 3D elastodynamics.

\subsection{Algorithms to perform   low-rank approximations}
Once the \emph{admissible} blocks are determined, an accurate rank-revealing algorithm is applied to determine low-rank approximations. Such an algorithm must be accurate (i.e. its result, the computed numerical rank, must be as small as possible) to avoid unnecessary computational costs.
 The   truncated Singular Value Decomposition (SVD)~\cite{golub2012matrix}
 gives the best low-rank approximation (Eckart-Young theorem) for unitary invariant norms (e.g. Frobenius or spectral norm).
Thus it 
 produces an approximation with the smallest possible numerical rank for a given prescribed accuracy. 
  But the computation of the SVD is expensive, i.e. in the order of $O(\max(m,n) \times \min(m,n)^2)$ for an  $m \times n$ matrix, and in addition it requires the computation of all the entries of  ${\mathbb A}$. In the context of the $\mathcal{H}$-matrices, the use of the SVD would induce the undesired need to assemble the complete matrix.
 
  The adaptive cross approximation (ACA)~\cite{bebendorf2015wideband,bebendorf2003adaptive} offers an interesting alternative to the SVD since it produces a quasi-optimal low-rank approximation without requiring the assembly of the complete matrix.
 The starting point of the ACA is  that every matrix of  rank $r$ is the sum of $r$ matrices of rank $1$. 
 The ACA is thus a greedy algorithm that improves the accuracy of the approximation by adding iteratively rank-1 matrices. At iteration $k$, the matrix is split into the rank $k$ approximation ${\mathbb B}_k= \sum_{\ell=1}^k {\bf u}_{\ell} {\bf v}_{\ell}^{*}=\mathbb{U}_k \mathbb{V}_k^*$ with  $\mathbb{U}_k \in \mathbb{C}^{m \times k}$, $ \mathbb{V}_k \in \mathbb{C}^{n \times k}$, and the residual ${\mathbb R}_k= {\mathbb A} - \sum_{\ell=1}^k {\bf u}_{\ell} {\bf v}_{\ell}^*$; ${\mathbb A}={\mathbb B}_k+{\mathbb R}_k$. 
The information is shifted iteratively from the residual to the approximant.  A stopping criterion
 is used to determine the appropriate rank to achieve the required accuracy. A straightforward choice of 
   stopping criterion is    
 \begin{equation}
\mbox{until} \quad  ||  {\mathbb A} - {\mathbb B}_k||_F \le \varepsilon_{\footnotesize\operatorname{ACA}}    || {\mathbb A}||_F
\label{stop_full}
\end{equation}
where $\varepsilon_{\footnotesize\operatorname{ACA}}>0$ is a given parameter, and $||.||_F$ denotes the Frobenius norm.
In the following, we denote $r_{\operatorname{ACA}}$ the numerical rank obtained by the ACA for a required accuracy $\varepsilon_{\footnotesize\operatorname{ACA}}$.
 The complexity of this algorithm  to generate an approximation of rank $r_{\operatorname{ACA}}$ is    $O(r_{\operatorname{ACA}}mn)$. 
 
 There are various ACAs that differ by the choice of the best pivot at each iteration. The simplest approach is the so-called
 fully-pivoted ACA and it 
 consists in choosing the pivot as the largest entry in the residual. But   similarly to the SVD, it requires the computation of all the entries of ${\mathbb A}$ to compute the pivot indices. It is not an interesting option for the construction of $\mathcal{H}$-matrices.
 
 The  partially-pivoted ACA proposes an alternative approach to choose the pivot avoiding the assembly of the complete matrix. The idea is to  maximize alternately  the residual  for only one of the two
indices  and to keep the other one  fixed. With this strategy,
 only one row and one column is  assembled at each iteration. More precisely,  
at iteration $k$, given ${\mathbb B}_k$ and assuming the row index $\hat{i}$ is known the algorithm is given by the following six steps:
 \begin{enumerate}
 \item Generation of the rows ${\bf a}:=  {\mathbb A}^* {\bf e}_{\hat{i}}$ and  $\begin{small}{\mathbb R}_k^* {\bf e}_{\hat{i}}={\bf a}-   \sum_{\ell=1}^k ({ u}_{\ell})_{\hat{i}}{\bf v}_{\ell}\end{small}$ 
 \item Find the  column index $\hat{j}:=\mbox{argmax}_{j}|({\mathbb R}_k)_{\hat{i}j}|$ and compute $\gamma_{k+1}= ({\mathbb R}_k)_{\hat{i},\hat{j}}^{-1}$
 \item Generation of the columns: ${\bf a}:=  {\mathbb A} {\bf e}_{\hat{j}}$  and ${\mathbb R}_k {\bf e}_{\hat{j}}={\bf a}-   \sum_{\ell=1}^k {\bf u}_{\ell} ({ v}_{\ell})_{\hat{j}}$ 
  \item Find the  next row index $\hat{i}:=\mbox{argmax}_{i}|({\mathbb R}_k)_{i\hat{j}}|$ 
  \item Compute vectors ${\bf u}_{k+1}:=\gamma_{k+1} {\mathbb R}_k {\bf e}_{\hat{j}}$, ${\bf v}_{k+1}:=  {\mathbb R}_k^* {\bf e}_{\hat{i}}$
  \item Update the approximation ${\mathbb B}_{k+1}={\mathbb B}_{k} + {\bf u}_{k+1} {\bf v}_{k+1}^*$
  \end{enumerate}
 
With this approach however, approximates and residuals are not computed explicitly nor stored, so the stopping criteria~(\ref{stop_full}) needs to be adapted. The common choice is a
 stagnation-based error estimate which is computationally inexpensive. The algorithm stops when the new rank-1 approximation does not improve the accuracy of the approximation. Since at each iteration $k$, ${\mathbb B}_k-{\mathbb B}_{k-1}= {\bf u}_{k}{\bf v}_{k}^*$, the stopping criteria now reads:
 \begin{equation}
\mbox{until } ||{\bf u}_{{k}}||_2 ||{\bf v}_{{k}}||_2 \le\varepsilon_{\operatorname{ACA}} ||{\mathbb B}_{{k}} ||_F.
\label{stop_part}
\end{equation}
The complexity of   the partially-pivoted   ACA is reduced to $O(r_{\operatorname{ACA}}^2 (m+n))$.
 Since the partially-pivoted   ACA  is a heuristic method, there exist
counter-examples where ACA fails~\cite{Borm2003} and variants  have been designed to improve the robustness of the method.
{Nevertheless,  in all the numerical examples we have performed,  we do not need them}. It is worth mentioning that other approaches exist such as fast multipole expansions~\cite{rokhlin1985rapid,greengard1987fast}, panel clustering~\cite{sauter2000variable,hackbusch1989fast}, quadrature formulas~\cite{borm2016approximation} or interpolations~\cite{messner2012fast}. These approaches combine approximation techniques and linear algebra. The advantages of the ACA are to be purely algebraic and easy to implement.

From now on, we shall write the result $\mathbb{B}_{r_{\operatorname{ACA}}}$ as $\mathbb{B}$ for short, and $\mathbb{B}\approx \mathbb{A}$ or  $\mathbb{A}\approx \mathbb{B}$. 

\subsection{Extension to problems with vector unknowns \label{ACA_vect}}
One specificity of this work is to consider $\mathcal{H}$-matrices and ACA in the context of systems of partial differential equations with vector unknowns.
There exist a lot of works both theoretical and numerical on the ACA for scalar problems,
in particular on the selection of non-zero pivots (since they are  used to normalize the new rank-1 approximation). 
Indeed for scalar problems, it is straightforward to find the largest non-zero entry in a given column. 
For      problems with vector unknowns in $\mathbb{R}^d$,  the system has a block structure, i.e.  each pair of nodes on the mesh does not define
a single entry  but rather a $d \times d$ subblock in the complete matrix. This happens for example for 3D elastodynamics where the Green's tensor is a $3 \times 3$ subblock.

Different strategies can be applied to perform the ACA on matrices with a block structure. The first strategy consists in ordering the system matrix such that it is composed of $9=3^2  $ subblocks of size $N_c \times N_c$ where $N_c$ is the number of points in the cloud. In 3D elastodynamics, this corresponds to the following partitioning of the matrix (below the solution is represented with ${\bf U}_{\omega}$):
\[
\left[
\begin{array}{ccc}
\mathbb{A}_{11} & \mathbb{A}_{12} & \mathbb{A}_{13}\\
\mathbb{A}_{21} & \mathbb{A}_{22} & \mathbb{A}_{23}\\
\mathbb{A}_{31} & \mathbb{A}_{32} & \mathbb{A}_{33}\\
\end{array}
\right],   \quad    (\mathbb{A}_{\alpha \beta})_{ij} = ({\bf U}_{\omega})_{\alpha \beta}(\x_i, \y_j) \ \ 1 \le \alpha,\beta \le 3, \ 1 \le i,j \le N_c.
\]
Then each submatrix is approximated independently with the conventional scalar ACA. This strategy is used 
in~\cite{MESSNER2010}.
It is well suited for iterative solvers.   But it cannot be adapted straightforwardly in the context of direct solvers since the recursive $2 \times 2$ block structure inherited from the binary tree is used. Indeed even though it is possible to determine the LU decomposition of each block of $\mathbb{A}$ (see Section~\ref{lu_solve}), such that
\[
\left[
\begin{array}{ccc}
\mathbb{A}_{11} & \mathbb{A}_{12} & \mathbb{A}_{13}\\
\mathbb{A}_{21} & \mathbb{A}_{22} & \mathbb{A}_{23}\\
\mathbb{A}_{31} & \mathbb{A}_{32} & \mathbb{A}_{33}\\
\end{array}
\right] \approx
\left[
\begin{array}{ccc}
\mathbb{L}_{11}\mathbb{U}_{11} & \mathbb{L}_{12}  \mathbb{U}_{12}& \mathbb{L}_{13} \mathbb{U}_{13}\\
\mathbb{L}_{21}\mathbb{U}_{21} & \mathbb{L}_{22} \mathbb{U}_{22}& \mathbb{L}_{23} \mathbb{U}_{23}\\
\mathbb{L}_{31} \mathbb{U}_{31}& \mathbb{L}_{32}  \mathbb{U}_{32} & \mathbb{L}_{33} \mathbb{U}_{33}\\
\end{array}
\right],
\]
on the other hand it would be very expensive to deduce the LU decomposition of $\mathbb{A}$ from these LU decompositions since what we are looking for  is a decomposition of the kind
\[
\left[
\begin{array}{ccc}
\mathbb{L}_{1}  & \mathbb{O} & \mathbb{O}\\
\mathbb{L}_{2} & \mathbb{L}_{3}  & \mathbb{O}\\
\mathbb{L}_{4} & \mathbb{L}_{5}  & \mathbb{L}_{6}  \\
\end{array}
\right]
\left[
\begin{array}{ccc}
 \mathbb{U}_{1} &  \mathbb{U}_{2}&  \mathbb{U}_{3}\\
\mathbb{O}  &  \mathbb{U}_{4}&   \mathbb{U}_{5}\\
\mathbb{O} & \mathbb{O}  &  \mathbb{U}_{6}\\
\end{array}
\right],\]
i.e. 12 factors instead of 18.
A solution would be to replace the binary by a ternary tree but in that case the clustering process would be more complex. 

The second naive approach consists in considering the complete matrix as a scalar matrix, i.e. to forget about the block structure. While appealing this approach fails in practice for 3D elastodynamics, whatever ordering is used, due to the particular structure of the matrix. Rewriting the Green's tensors~(\ref{elasto_U})-(\ref{elasto_T}) component-wise~\cite{bonnet1999boundary}  and using the Einstein summation convention, it reads
\begin{multline*}
\hspace*{-0.5cm} {({\bf U}_{\omega})_{\alpha \beta}({\x},{\y})=\frac{1}{4 \pi  \varrho} (a_1 \delta_{\alpha \beta}+a_2 \varrho_{,\alpha}  {\varrho}_{,\beta}),} \\
 { ({\bf T}_{\omega})_{\alpha \beta}({\x},{\y})=\frac{1}{4 \pi \varrho} \Big[2a_3  {\varrho}_{,\alpha}  {\varrho}_{,\gamma} {\varrho}_{,\beta} +a_4(\delta_{\alpha \beta} \varrho_{,\gamma}
+ \delta_{\gamma \beta} \varrho_{,\alpha}) +a_5 \delta_{\alpha \gamma}\varrho_{,\beta} \Big] { n}_{\gamma}(\y)}
\end{multline*}
where the constants $a_1, ..., \ a_5$ depend only on the mechanical properties, $ {\varrho}=||{\x}-{\y}||$ and $ {\varrho}_{,\alpha}=\frac{\partial }{\partial y_\alpha}\varrho({\x},{\y})$.
As a result, as soon as ${\x}$ and ${\y}$ belong to the same plane (let's say $x_3=y_3$ for simplicity) the Green's tensors simplify to
\begin{equation}
\hspace*{-0.4cm}{\bf U}_{\omega}({\x},{\y}) = \frac{1}{4 \pi   {\varrho}}\left[
\begin{array}{ccc}
a_{1}+a_2 {\varrho}^2_{,1} & a_{2} {\varrho}_{,1} {\varrho}_{,2} & 0\\
a_{2}  {\varrho}_{,1} {\varrho}_{,2} & a_{1}+a_2  {\varrho}^2_{,1}  & 0\\
0 & 0 & a_{1}\\
\end{array}
\right], \ {\bf T}_{\omega}({\x},{\y}) = \frac{1}{4 \pi  {\varrho}}\left[
\begin{array}{ccc}
0 & 0& a_4  {\varrho}_{,1} \\
0 & 0& a_4  {\varrho}_{,2}\\
 a_4  {\varrho}_{,1} & a_4  {\varrho}_{,2}&0\\
\end{array}
\right].
\label{FSplane}
\end{equation}
It is then clear that the complete matrix ${\mathbb A}$ composed of such subblocks is a reducible matrix (there are decoupled groups of unknowns after discretization). While the fully-pivoted ACA will succeed in finding the best pivot to perform low-rank approximations, the partially pivoted ACA will only cover parts of the matrix resulting in a non accurate approximation of the initial matrix. As a result, this approach cannot be applied for 3D elastodynamics. 
Another strategy would be to adapt the choice of the pivot to reducible matrices, i.e. to cover all the decoupled groups of unknowns. This is possible with the use of the ACA+ algorithm that defines a reference column and a reference row along whose the pivots are looked for~\cite{borm2003introduction}. Even though this approach could fix the problem it does not seem to be the best suited one since it does not use the vector structure of 3D elastodynamic problems.

To sum up, it is now important to use an algorithm that takes into account the particular structure of our matrix such that the vector ACA does not rely anymore on a rank-1 update but instead on a rank-3 update. The central question is then how to find the pivot used for this rank-3 update instead. There are 3 possible strategies:
\begin{enumerate}
\item To look for the largest scalar pivot, determine the corresponding point in the cloud and update all  3 DOFs linked to this point simultaneously. This approach is not stable since the $3 \times 3$ subblock pivot may not be invertible  (cf. a discretized version of~(\ref{FSplane})).
\item The second strategy is to look for the $3\times 3$ subblock with the largest norm. Again this approach fails in practice since  the $3 \times 3$ subblock pivot may not be invertible.
\item The third strategy used in this work is to compute the singular values $\sigma_1\ge \sigma_2 \ge \sigma_3$ of every candidate subblock. In order to achieve convergence and to avoid singular pivots, the safest approach consists in choosing the pivot corresponding to the subblock with the largest $\sigma_3$. A similar approach is used for electromagnetism in~\cite{rjasanow2016matrix}.
\end{enumerate}

\begin{remark}
It is worth noting that for some specific configurations the 3D elastodynamic double layer potential ${\bf T}_{\omega}$ may lead to a  matrix with only singular subblocks. In such cases, the randomized SVD~\cite{liberty2007randomized} is preferred to the vector ACA.
\end{remark}

    \section{$\mathcal{H}$-matrix based iterative and direct solvers}
    \subsection{$\mathcal{H}$-matrix based iterative solver \label{iter_solver}}
 Once the $\mathcal{H}$-matrix representation of a matrix is computed, it is easy to derive an $\mathcal{H}$-matrix based iterative solver. The only operation required
 is an efficient  matrix-vector product. It is performed hierarchically by going through the block cluster tree $\mathcal{P}$. At the leaf level, there are two possibilities. If
 the block of size $m \times n$ does not admit a low-rank approximation (\emph{non-admissible} block), then the standard matrix-vector product is used with cost $O(mn)$. Otherwise, 
  the block is marked as \emph{admissible} such that   a low-rank approximation has been computed~: ${\mathbb A}_{\tau \times \sigma} \approx \mathbb{B}_{\tau \times \sigma}$. The cost of this part of the matrix-vector product is then reduced from $O(mn)$   to $O(r_{\operatorname{ACA}}(m+n))$ where $r_{\operatorname{ACA}}$ is the  numerical rank of the block $\mathbb{B}_{\tau \times \sigma}$ computed with the ACA.

  \subsection{$\mathcal{H}$-LU factorization and direct solver \label{lu_solve}}
One of the advantages of the $\mathcal{H}$-matrix representation is the possibility to derive a fast direct solver. Due to the hierarchical block structure of a  $\mathcal{H}$-matrix, the LU factorization is performed recursively on  $2 \hspace*{-0.1cm} \times  \hspace*{-0.1cm}2$    block matrices of the form
   \begin{displaymath}
	\left( 
   	\begin{array}{ c  c } 
                   {\mathbb A}_{\tau_1 \times \sigma_1} &  {\mathbb  A}_{\tau_1 \times \sigma_2}\\              
                    {\mathbb  A}_{\tau_2 \times \sigma_1} &  {\mathbb  A}_{\tau_2 \times \sigma_2}
    	\end{array} 
	\right)=\left( 
   	\begin{array}{ c  c } 
                   {\mathbb  L}_{\tau_1 \times \sigma_1} & \mathbb{O}\\ 
                  {\mathbb L}_{\tau_2 \times \sigma_1} & {\mathbb  L}_{\tau_2 \times \sigma_2}
    	\end{array} 
	\right) \left( 
   	\begin{array}{ c  c } 
                   {\mathbb U}_{\tau_1 \times \sigma_1} &  {\mathbb U}_{\tau_1 \times \sigma_2}\\ 
                    \mathbb{O} &  {\mathbb U}_{\tau_2 \times \sigma_2}
    	\end{array} 
	\right) 
	\end{displaymath}
For such block matrices, we recall that the LU factorization is classically decomposed into 4 steps:
\begin{enumerate}
\item LU decomposition to compute  ${\mathbb  L}_{\tau_1 \times \sigma_1}$   and  ${\mathbb U}_{\tau_1 \times \sigma_1}$:  ${\mathbb A}_{\tau_1 \times \sigma_1} =   {\mathbb  L}_{\tau_1 \times \sigma_1}     {\mathbb U}_{\tau_1 \times \sigma_1}$;
\item Compute $  {\mathbb U}_{\tau_1 \times \sigma_2}$ from    $  {\mathbb  A}_{\tau_1 \times \sigma_2}={\mathbb  L}_{\tau_1 \times \sigma_1} {\mathbb U}_{\tau_1 \times \sigma_2}$;  
\item Compute $   {\mathbb L}_{\tau_2 \times \sigma_1} $ from $    {\mathbb  A}_{\tau_2 \times \sigma_1}=     {\mathbb L}_{\tau_2 \times \sigma_1} {\mathbb U}_{\tau_1 \times \sigma_1}$;
\item LU decomposition to compute ${\mathbb  L}_{\tau_2 \times \sigma_2}$   and  ${\mathbb U}_{\tau_2 \times \sigma_2}$:  $   {\mathbb  A}_{\tau_2 \times \sigma_2}-   {\mathbb L}_{\tau_2 \times \sigma_1} {\mathbb U}_{\tau_1 \times \sigma_2}= {\mathbb  L}_{\tau_2 \times \sigma_2}  {\mathbb U}_{\tau_2 \times \sigma_2}$. 
 \end{enumerate}
The obvious difficulty is the fact that the process is recursive (cf steps~1 and~4). In addition, Step~4 requires the addition and multiplication of blocks which may come in different formats. In more general terms, the difficulty is to evaluate
\[
 {\mathbb A}^{'''}_{\tau \times \sigma} \leftarrow {\mathbb A}_{\tau \times \sigma} +    {\mathbb A}^{'}_{\tau \times \sigma'}  {\mathbb A}^{''}_{\tau' \times \sigma}.
\]
 Since each block may belong to three different classes this results into $27$ cases to take into account.
 For example, if  ${\mathbb A}_{\tau \times \sigma} $ is a full matrix,  $ {\mathbb A}^{'}_{\tau \times \sigma'} $ an $\mathcal{H}$-matrix, i.e. a \emph{non-admissible} block at a non-leaf level  and ${\mathbb A}^{''}_{\tau' \times \sigma}$ a matrix having a low-rank representation, then the computation is efficiently decomposed into the following steps. First, 
  ${\mathbb A}^{''}_{\tau' \times \sigma} $ being a low-rank matrix  it is decomposed into  ${\mathbb A}^{''}_{\tau' \times \sigma} \approx \mathbb{B}^{''}_{\tau' \times \sigma}= \mathbb{U} \mathbb{V}^*$ and 
using the fast $\mathcal{H}$-matrix / vector product   the product between ${\mathbb A}^{'}_{\tau \times \sigma'}$ and each 
column  of  ${\mathbb U}$ is evaluated. In a second step, the result is multiplied to ${\mathbb V}^{*}$ and the final matrix is directly a full matrix.
 A more detailed review of the 27 cases is presented in~\cite{desiderio}.

It is known that the numerical rank obtained by the ACA is not optimal.
Similarly the numerical rank for blocks obtained after an  addition and a multiplication and stored in a low-rank or $\mathcal{H}$-matrix format  is not optimal.
All these blocks are further  recompressed to optimize the numerical ranks.
 Starting from the low-rank approximation $\mathbb{B}_{\tau \times \sigma}={\mathbb U} {\mathbb V}^{*} $,  the reduced QR-decompositions of ${\mathbb U}$ and $ {\mathbb V}$ are computed, i.e.
 \[
\mathbb{B}_{\tau \times \sigma} = {\mathbb U} {\mathbb V}^{*} = {\mathbb Q}_U {\mathbb R}_U {\mathbb R}_V^* {\mathbb Q}_V^*.
 \]
Recall that the only decomposition that yields the optimal numerical rank is   the truncated singular value decomposition. It is thus applied to  the reduced matrix ${\mathbb R}_U {\mathbb R}_V^* \approx \mathbb{PSL}^*$
 such that the cost of this step is reasonable.
  \[
\mathbb{B}_{\tau \times \sigma} = {\mathbb Q}_U {\mathbb R}_U {\mathbb R}_V^* {\mathbb Q}_V^*= ({\mathbb Q}_U {\mathbb P} {\mathbb S}^{1/2})   ({\mathbb S}^{1/2}  {\mathbb L}^* {\mathbb Q}_V^*)= {\mathbb U'} {\mathbb V'}^{*}.
 \]
During the construction of the $\mathcal{H}$-matrix representation, the   parameter $\varepsilon_{\operatorname{ACA}}$ is also used to determine the number of singular values selected during the recompression process. Moreover, in the case of a recompression after an addition and a multiplication, i.e. when the format of a block is modified, we introduce the parameter $\varepsilon_{\operatorname{LU}}$ to determine the level of accuracy required. 

\subsection{Proposition of an   estimator to certify the results of the $\mathcal{H}$-LU direct solver \label{esti_direct}}

Since the $\mathcal{H}$-LU direct solver is based on a heuristic method (the ACA) to perform the low-rank approximations and due to the various modifications on the  $\mathcal{H}$-matrix performed during the LU factorization, it is relevant to propose a simple and efficient way to certify the  obtained results. To this aim, we propose an   estimator. We consider the initial system
\[
\mathbb{A}{\bf x}={\bf b}.  
\]
 We denote ${\mathbb A}_{\mathcal{H}}$ the $\mathcal{H}$-matrix representation of ${\mathbb A}$ (in which low-rank approximations have been performed) and  ${\mathbb L}_{\mathcal{H}}{\mathbb U}_{\mathcal{H}}\approx {\mathbb A}_{\mathcal{H}}$ its  $\mathcal{H}$-LU factorization. The solution of the approximated system is ${\bf x}_0$~: ${\mathbb L}_{\mathcal{H}}{\mathbb U}_{\mathcal{H}} {\bf x}_0={\bf b}$. Our aim is to give an upper bound of $||{\bf b}- {\mathbb A}{\bf x}_0||_2$.
 Since
 \[
 {\bf b}-{\mathbb A}{\bf x}_0={\bf b}-{\mathbb A}_{\mathcal{H}}{\bf x}_0+{\mathbb A}_{\mathcal{H}}{\bf x}_0-{\mathbb A}{\bf x}_0
 \]
this yields for $\alpha\in \{2,F\}$\footnote{Recall that $||{\mathbb B}||_2 \le ||{\mathbb B}||_F$ for all square matrices.} 
\begin{equation}
\begin{array}{rl}
 \displaystyle \frac{||{\bf b}- {\mathbb A}{\bf x}_0||_2}{||{\bf b}||_2}& \displaystyle\le  \frac{1}{||{\bf b}||_2} (||{\bf b}- {\mathbb A}_{\mathcal{H}} {{\bf x}_0}||_2+|| {\mathbb A}_{\mathcal{H}}-{\mathbb A}||_{\alpha} ||{{\bf x}_{0}}||_2) { =   \frac{1}{||{\bf b}||_2}  (\delta+\delta_{\mathcal{H},\alpha}  ||{\bf x}_{0}||_2)}\\[2ex]
& \mbox{where } \delta=||{\bf b}- {\mathbb A}_{\mathcal{H}}{{\bf x}_0}||_2 \mbox{ and } \delta_{\mathcal{H},\alpha}=||{\mathbb A}_{\mathcal{H}}-{\mathbb A}||_{\alpha}.
\end{array}
\label{est_direct}
\end{equation}
 It is worth noting that $\delta_{\mathcal{H},\alpha}=||{\mathbb A}_{\mathcal{H}}-{\mathbb A}||_{\alpha} $ estimates the accuracy of the  $\mathcal{H}$-matrix representation (i.e. the influence of the parameter $\varepsilon_{\operatorname{ACA}}$ in the computation of the low-rank approximations) while 
 $\delta=||{\bf b}- {\mathbb A}_{\mathcal{H}}{{\bf x}_0}||_2$ accounts for the stability of the $\mathcal{H}$-LU factorization (i.e. the influence of the parameter $\varepsilon_{\operatorname{LU}}$ in the LU factorization). The evaluations of $\delta$, $||{\bf b}||_2$ and $||{\bf x}_0||_2$ reduce to the computation of the norm of a vector.  {The most expensive part is $\delta_{\mathcal{H},\alpha}$ but its
 evaluation is performed with the Frobenius norm to reduce the cost.} Also  this term does not depend on the right hand side meaning that for multiple right hand sides this error estimate is not expensive.

\section{3D Elastodynamics Boundary Element Method}

Let's consider a bounded domain $\Omega^-$ in $\R^3$ representing an obstacle with a closed Lipschitz
boundary $\Gamma:=\partial \Omega^-$. Let  $\Omega^+$ denote the associated exterior domain $\R^3\backslash\overline{\Omega^-}$ and $\nn$ its  outward
unit normal vector field on its  boundary $\Gamma$.  We consider  
the propagation of time-harmonic 
waves in a three-dimensional isotropic and homogeneous elastic medium  modeled by the Navier equation~(\ref{NE}).
The following results about traces of vector fields  and integral
representations of time-harmonic elastic fields can be found in~\cite{darbas2015well}.

\subsection{Traces, integral representation formula and integral equations}
We denote by $H^s_{loc}({\Omega^\pm})$ and $H^s(\Gamma)$ the standard (local in the case of the exterior domain) complex valued Hilbertian Sobolev spaces of order $s\in\R$ ($|s|\le 1$ for  $H^s(\Gamma)$) defined on ${\Omega^\pm}$ and $\Gamma$ respectively (with the convention $H^0=L^2$).  Spaces of vector functions will be denoted by boldface letters, thus $\HH^s=(H^s)^3$.  We  set $\mathbf{\mathbf{\Updelta}}^*\uu:=\Div\upsigma(\uu)=(\lambda+2\mu)\vnabla\Div\uu-\mu\Rot\Rot\uu$
 and introduce the energy spaces
$\HH^1_{+}(\mathbf{\Updelta}^*) := \big\{\uu\in \HH^1_{loc}({\Omega^+}):\;\mathbf{\Updelta}^*\uu\in \LL_{loc}^2({\Omega^+})\big\}$ and
$\HH^1_{-}(\mathbf{\Updelta}^*) := \big\{\uu\in \HH^1({\Omega^-}):\;\mathbf{\Updelta}^*\uu\in \LL^2({\Omega^-})\big\}$. The traction trace for elastodynamic problems is defined
by $\tt_{|\Gamma}:=\TT\uu$ where $\TT$ is the traction operator defined in~(\ref{elasto_T}). We recall that we have $\uu_{|\Gamma}\in\HH^{\frac{1}{2}}(\Gamma)$ and $\tt_{|\Gamma}\in\HH^{-\frac{1}{2}}(\Gamma)$  for all $\uu\in\HH^1_{\pm}(\mathbf{\Updelta}^*)$.

For a solution $\uu\in\HH^1_{+}(\mathbf{\Updelta}^*)$ to the Navier equation \eqref{NE} in $\Omega^+$, that satisfies the Kupradze radiation conditions, the 
Somigliana integral representation of the field is given by 
\begin{equation}
\label{somigliana}
\uu(\x)=\,\mathcal{D} \uu_{|\Gamma}(\x)\,-\,\mathcal{S}\tt_{|\Gamma}(\x), \quad \x \in \Omega^+,
\end{equation} 
where the single- and  double-layer potential operators are respectively defined by 
\begin{equation}
\label{potential}
\mathcal{S} \vphi(\x)=\int_{\Gamma}{\bf U}_{\omega}(\x,\y) \vphi(\y)ds(\y)\;\;\text{ and }\;\; \mathcal{D} \vpsi(\x)=\int_{\Gamma}\transposee{\left[\TT_{\y}{\bf U}_{\omega}(\x,\y)\right]}\vpsi(\y)ds(\y) \quad \x \in \mathbb{R}^3 \backslash \Gamma.\end{equation}

The potentials $\mathcal{S}$ (resp. $\mathcal{D}$) are continuous from $\HH^{-1/2}(\Gamma)$  to $\HH^1_{-}(\mathbf{\Updelta}^*) \cup \HH^1_{+}(\mathbf{\Updelta}^*)$ (resp. from 
 $\HH^{1/2}(\Gamma)$    to $\HH^1_{-}(\mathbf{\Updelta}^*) \cup \HH^1_{+}(\mathbf{\Updelta}^*)$). For any $\vphi \in \HH^{-1/2}(\Gamma)$ and $\vpsi \in \HH^{1/2}(\Gamma)$, the potentials $\mathcal{S}\vphi$ and $\mathcal{D}\vpsi$
solve the Navier equation in $\Omega^+$ and $\Omega^-$, and satisfy the Kupradze radiation condition.
The exterior and interior Dirichlet $(\gamma^{\pm}_0)$ and traction  $(\gamma^{\pm}_1)$ traces of $\mathcal{S}$ and $\mathcal{D}$ are given by
 \[
 \gamma^{\pm}_0 \mathcal{S} =S, \quad  \gamma^{\pm}_1 \mathcal{S} =\mp \frac{1}{2}I + D', \quad  \gamma^{\pm}_0 \mathcal{D} =\pm \frac{1}{2}I + D
 \]
 where the operators $S$ (resp. $D$) are continuous from $\HH^{-1/2}(\Gamma)$ to $\HH^{1/2}(\Gamma)$ (resp. continuous from  $\HH^{1/2}(\Gamma)$ to $\HH^{-1/2}(\Gamma)$) and are given by
 \begin{equation}
{S} \vphi(\x)=\int_{\Gamma}{\bf U}_{\omega}(\x\,,\y) \vphi(\y)ds(\y)\;\;\text{ and }\;\; {D} \vpsi(\x)=\int_{\Gamma}\transposee{\left[\TT_{\y}{\bf U}_{\omega}(\x\,,\y)\right]}\vpsi(\y)ds(\y), \ \x \in \Gamma.\end{equation}
 The scattering problem is formulated as follows : Given an incident wave field $\uinc$ which is assumed to solve the Navier equation in the absence of any scatterer, find the displacement $\uu$ solution to the Navier equation \eqref{NE} in $\Omega^+$ which satisfies the  Dirichlet boundary condition on $\Gamma$
\begin{equation}\label{Dirichlet}\uu_{\vert\Gamma}+ \uu^{\mathrm{inc}}_{\vert\Gamma} = 0.\end{equation}
Applying the potential theory, the elastic scattering problem reduces to a boundary integral equation
 \begin{equation}
 \label{EFIE}
S  (\tt_{| \Gamma}+\tt^{inc}_{| \Gamma}) (\x)= 
 \uu^{inc}_{| \Gamma}({\x}) ,  \quad {\x} \in \Gamma 
\end{equation}
 \begin{equation}
 \label{MFIE}
\mbox{or} \quad  (\frac{I}{2} + D')  (\tt_{| \Gamma}+\tt^{inc}_{| \Gamma}) (\x)= 
 \tt^{inc}_{| \Gamma}({\x}) ,  \quad {\x} \in \Gamma.  
\end{equation}
In the following, $\mathcal{H}$-matrix based solvers are applied and studied in the special case of  elastodynamic scattering problems
but the method can  be applied to the solution of any boundary integral equation defined in terms of the single and double layer potential operators.

\subsection{Classical concepts of the Boundary Element Method}
The main ingredients of the Boundary Element Method are a transposition of  the concepts developed for  the Finite Element Method~\cite{bonnet1999boundary}.
First, the numerical solution of the boundary integral equations~(\ref{EFIE}) or (\ref{MFIE}) is
based on a discretization of the surface $\Gamma$ into
$N_{E}$ isoparametric boundary elements of order one, i.e.  {three-node triangular elements}. 
Each physical element $E_e$ on the approximate boundary is mapped onto a reference element $\Delta_e$  via an affine mapping
\[
{\boldsymbol \xi} \in \Delta_e \rightarrow \y({\boldsymbol \xi})  \in E_e, \quad 1 \le e \le N_e.
\]
$ \Delta_e$ is the reference triangle in the $(\xi_1,\xi_2)$-plane. The $N_c$ interpolation points $\y_1, \ldots, \y_{N_c}$ are chosen as the vertices of the mesh.
Each component of the total traction field is  approximated with globally  continuous, piecewise-linear shape functions $(v_i(\y))_{1\le i\le N_c}$: $v_i(\y_j)=\delta_{ij}$ for $1 \le i,j \le N_c$.  A boundary element $E_e$ contains exactly $3$ interpolation nodes $(\y_k^{e})_{1 \le k \le 3}$ associated with $3$ basis functions $(v^e_k)_{1 \le k \le 3}$. 
These basis functions are related to the canonical basis $(\hat{v}_k)_{1 \le k \le 3}$ defined on the reference element $\Delta_e$ by $v_k^e(\y({\boldsymbol \xi}))=\hat{v}_k({\boldsymbol \xi})$.
Each component of the total traction field ${\boldsymbol p}(\y) =  (\tt_{| \Gamma}+\tt^{inc}_{| \Gamma}) (\y)$ is approximated on the element $E_e$ by
\[
{\boldsymbol p}_{\alpha}(\y) \approx \sum_{k=1}^{3} p^{k}_{\alpha} v_k^e({\y}) \quad (1 \le \alpha \le 3),
\]
where   $p^k_{\alpha}$ denotes  the approximation of the nodal value of the component $\alpha$ of the vector ${\boldsymbol p}(\y_k)$.
To discretize the boundary integral equations~(\ref{EFIE}) or (\ref{MFIE}) we consider the collocation approach. It consists in enforcing the equation at a finite number of collocation points $\x$.  To have a solvable discrete problem, one has to choose $N_c$ collocation points. The $N_{c}$ traction
approximation nodes thus defined also serve as collocation points, i.e. $(\x_i)_i=(\y_j)_j$.
This discretization process transforms~(\ref{EFIE}) or (\ref{MFIE}) into a square
complex-valued  linear system of size $3N_{c}$ of the form
\begin{equation} {\mathbb A} {\bf p} = {\bf b}, \label{discre_syst}
\end{equation}
where the $(3N_c)$-vector ${\bf p}$ collects the   degrees of freedom
(DOFs), namely the nodal traction components, while  the $(3N_c)$-vector ${\bf b}$
arises from the imposed incident wave field. Assembling
the matrix ${\mathbb A}$ classically~\cite{bonnet1999boundary} requires the computation of all element
integrals for each collocation point, thus requiring a computational
time of order $O(N_c^{2})$.

It is interesting to note that   a data-sparse representation 
of the matrix ${\mathbb A}$  is a direct consequence of the discretization of the Green's tensor.
For example with the defined process, the discretization of equation~(\ref{EFIE})   leads to the system
  \[
 \hspace*{-1cm}
{\mathbb U}_{\omega}
 {\mathbb W}_y
 {\mathbb V}_y {\bf p}={\bf b} 
\]
%
where the \emph{diagonal} matrix ${\mathbb W}_y$ corresponds to the weights used to evaluate numerically  the integrals, the matrix ${\mathbb V}_y$ corresponds to shape functions evaluated at the quadrature points of the reference element and the matrix  ${\mathbb U}_{\omega}$ corresponds to the evaluation of the Green's tensor at the collocation points $(\x_i)_{i=1,\ldots,N_c}$ and interpolation points $(\y_j)_{j=1,\ldots,N_c}$.
 The relevant question is then how the $\mathcal{H}$-matrix representation of ${\mathbb U}_{\omega}$ can be transmitted to $\mathbb{A}$.
 The key point is to remark that the matrix  ${\mathbb V}_y$ is a sparse matrix whose connectivity can be described as follows: it has a non-zero entry if there exists one triangle that has the two corresponding nodes as vertices.
 From that argument, it appears that if no efficient $\mathcal{H}$-matrix representation of $\mathbb{U}_{\omega}$ is available, we cannot expect to find an 
  efficient $\mathcal{H}$-matrix representation of $\mathbb{A}$. On the other hand, since the structure of the matrix ${\mathbb V}_y$ is based on a notion of distance similarly to the construction of the binary tree $\mathcal{T}_I$, we can expect that if an efficient $\mathcal{H}$-matrix representation of $\mathbb{U}_{\omega}$ is available, an efficient 
   $\mathcal{H}$-matrix representation of $\mathbb{A}$ will be found. However,   the interface separating two subdomains used to compute the binary tree can also separate nodes belonging to the same elements. It is thus very likely that the ranks of the subblocks of $\mathbb{A}$ will be larger than the ranks of the same blocks in $\mathbb{U}_{\omega}$.
 Importantly,  this property is independent of the order of approximation used in the BEM.

As a consequence, although the behavior of the $\mathcal{H}$-matrix based solvers for 3D elastodynamics are presented in the context of the BEM, the observations can  nevertheless be extended to other configurations where a Green's tensor is discretized over a cloud of points.

\section{Application of $\mathcal{H}$-matrices to oscillatory kernels: theoretical estimates \label{theory}}


The efficiency of $\mathcal{H}$-matrix based solvers depends on the possible storage reduction obtained by low-rank approximations.
It is thus important to estimate the memory requirements of the method in the context of 3D elastodynamics. We follow the proof of~\cite{hackbusch2015hierarchical} proposed in the context of asymptotically smooth kernels.

\subsection{Storage estimate}
We use the notations introduced in Section~\ref{section_tree}.
Since we are using a balanced binary tree $\mathcal{T}_I$, the number of levels is $L(I)\le \log_2|I| - \log_2 N_{\operatorname{leaf}}+1$. We denote by  $N_c=|I|$  the number of points in the cloud.
For a \emph{non-admissible} block $\tau \times \sigma \in \mathcal{P}^{\operatorname{non-ad}}$, the memory requirements are 
\[
N_{st}= |\tau| \ |\sigma| = \mbox{min} \{|\tau|, |\sigma|  \} \ \mbox{max} \{|\tau|, |\sigma|  \} \mbox{ for } \tau \times \sigma \in \mathcal{P}^{\operatorname{non-ad}}.
\]
Since \emph{non-admissible} blocks can only appear  at the leaf level,  the memory requirements are bounded by
\[
N_{st}  \le  N_{\operatorname{leaf}} \ \max \{|\tau|, |\sigma|  \} \le N_{\operatorname{leaf}} (|\tau|+ |\sigma|)   \mbox{ for }  \tau \times \sigma \in \mathcal{P}^{\operatorname{non-ad}}.
\]
For \emph{admissible} blocks, the memory requirements are
\[
N_{st}= r_{\tau \times \sigma} (|\tau| + |\sigma|) \le  r_{\operatorname{ACA}}^{\max} (|\tau|+ |\sigma|) \mbox{ for }  \tau \times \sigma\in \mathcal{P}^{\operatorname{ad}},
\]
where  $r_{\tau \times \sigma}$ is the rank obtained by the ACA, resp.  $r_{\operatorname{ACA}}^{\max}$ denotes the maximum rank obtained by the ACA among all the \emph{admissible} blocks.
The total memory requirement is  given by
\[
N_{st}({\mathbb A})=\sum_{\tau \times \sigma \in \mathcal{P}^{\operatorname{ad}} } r_{\tau \times \sigma} (|\tau| + |\sigma|)  + \sum_{\tau \times \sigma \in \mathcal{P}^{\operatorname{non-ad}}  }  |\tau| \ |\sigma| \le \mbox{max} \{N_{leaf}, r_{ACA}^{\max}\}\sum_{\tau \times \sigma\in \mathcal{P} }   (|\tau|+ |\sigma|).
\]
It is clear that for a binary tree
\[
\sum_{\tau \in \mathcal{T}_i} |\tau| = \sum_{\ell=0}^{L(I)-1} |I^{(\ell)}|=L(I) |I|.
\]
Then, using the definition of the sparsity pattern we obtain
 \[
\sum_{\tau \times \sigma \in \mathcal{P}}|\tau| = \sum_{\tau \in \mathcal{T}_I} \Big[ |\tau| \sum_{\sigma | \tau \times \sigma \in \mathcal{P} } 1 \Big] \le C_{sp} \sum_{\tau \in \mathcal{T}_I} |\tau| \le  C_{sp} \ L(I) \ |I|.
\]
 A similar bound is found for $\sum_{\tau \times \sigma \in \mathcal{P}}|\sigma|  $ such that the  storage requirement is bounded by 
 \[
 N_{st}({\mathbb A}) \le 2 C_{sp} \max \{r_{\operatorname{ACA}}^{\max},N_{\operatorname{leaf}}\} N_c (\log_2 N_c - \log_2 N_{\operatorname{leaf}} +1).\]
In this storage estimate, the parameters $N_c$ and $ N_{\operatorname{leaf}}$ are known.
The sparsity pattern $C_{sp}$ depends on the partition $\mathcal{P}$ and will be discussed later (Section~\ref{impl_issues}).
In the context of oscillatory kernels, $r_{\operatorname{ACA}}^{\max}$ is known to depend on the frequency. The aim of the rest of this section is to characterize this dependence. We consider first   the 3D Helmholtz Green's function.

 \subsection{Convergence of the Taylor expansion of the  3D Helmholtz Green's function} 
 It is well-known that $\mathcal{H}$-matrices are efficient representations for matrices coming from the discretization of asymptotically smooth kernels,
  i.e.  kernels satisfying  Definition~\ref{asym_kern} below. For such matrices the method is well-documented and estimates  are provided for the Taylor expansion and interpolation errors.
  
  \begin{definition}
A kernel $s(.,.)~: \mathbb{R}^3 \times \mathbb{R}^3 \rightarrow \mathbb{R} $ is asymptotically smooth if there exist two constants $c_1,c_2$ and a singularity degree $\sigma \in \mathbb{N}_0$ 
 such that $\forall z \in \{x_{\alpha},y_{\alpha}\}$, $ \forall n \in \mathbb{N}_0$, $\forall \x \neq \y$
   \[
 |\partial_z^n s(\x,\y)| \le n! c_1 (c_2 ||\x-\y||)^{-n -\sigma}.
 \]
 \label{asym_kern}
 \end{definition}
 The capability  to produce low rank approximations is closely related to the concept of  degenerate functions~\cite{bebendorf2008hierarchical}.
A kernel function is said to be degenerate   if it can be well approximated (under some assumptions) by a  sum of functions with separated variables.
In other words, noting $X$ and $Y$ two domains of $\mathbb{R}^3$,   we are looking for an approximation $s^r $  of  $s$ on $X \times Y$ with $r$ terms  such that it writes
\[
s^r(\x,\y)=\sum_{\ell=1}^{r} u^{(r)}_{\ell}(\x) v^{(r)}_{\ell}(\y) \quad \x \in X, \ \ \y \in Y.
\]
$r$ is called the \emph{separation rank} of $s^r $ and $s(\x,\y)=s^r(\x,\y)+R_s^r(\x,\y)$ with $R^r_s$ the remainder. 
Such an approximation is obtained for instance with 
a Taylor expansion with $r:=| \{{\boldsymbol \alpha}\in \mathbb{N}^3_0:  |{\boldsymbol \alpha}| \le m\}|$ terms, i.e. 
\[
s^r(\x,\y)=\sum_{{\boldsymbol \alpha}\in \mathbb{N}^3_0, |{\boldsymbol \alpha}| \le m }   (\x-\x_0)^{\boldsymbol \alpha} \frac{1}{{\boldsymbol \alpha}!} \partial_x^{\boldsymbol \alpha}s(\x_0,\y) + R^r_s \quad \x \in X, \ \ \y \in Y
\]
where $\x_0\in X$ is the centre of the expansion.

A kernel function is said to be a separable expression in $X \times Y$ if the remainder converges quickly to zero. The main result for asymptotically smooth kernels is  presented in~\cite[Lemma 3.15]{bebendorf2008hierarchical} and \cite[Theorem 4.17]{hackbusch2015hierarchical}. For such kernels, it can be easily shown that
\[
|R^r_s (\x,\y)| \le C' \sum_{\ell=m}^{\infty} \Big( \frac{ ||\x-\x_0||}{c_2 ||\x_0-\y||}\Big)^{\ell}
\]
where $c_2$ is the constant appearing in Definition~\ref{asym_kern}. 
The convergence of $R^r_s$ is thus ensured  for all $\y \in Y$ such that
\begin{equation}
\gamma_x:= \frac{\max_{\x \in X} ||\x-\x_0||}{c_2 ||\x_0-\y||}< 1.
\label{admi_cond_taylor_h}
\end{equation}
Provided that  the condition (\ref{admi_cond_taylor_h}) holds, the remainder is bounded by
\begin{equation}
\label{error_est_remain}
\begin{array}{rl}
|R^r_s(\x,\y)|\le &\displaystyle{C'  
 \frac{\gamma_x^m}{1-\gamma_x}  \xrightarrow[m \to \infty]{}  0.}
 \end{array}
\end{equation}
In this configuration, the Taylor expansion of the asymptotically smooth  kernel converges exponentially with convergence rate 
$\gamma_x$. Since {$r$ is the cardinal of the set $\{{\boldsymbol \alpha}\in \mathbb{N}^3_0:  |{\boldsymbol \alpha}| \le m\}$, it holds } $m \sim r^{1/3}$, and then
$r\approx |\log \varepsilon|^3$ is expected to achieve an approximation with accuracy $\varepsilon>0$.

\begin{figure}
\begin{center}
\begin{tabular}{c}
\includegraphics[scale=0.35]{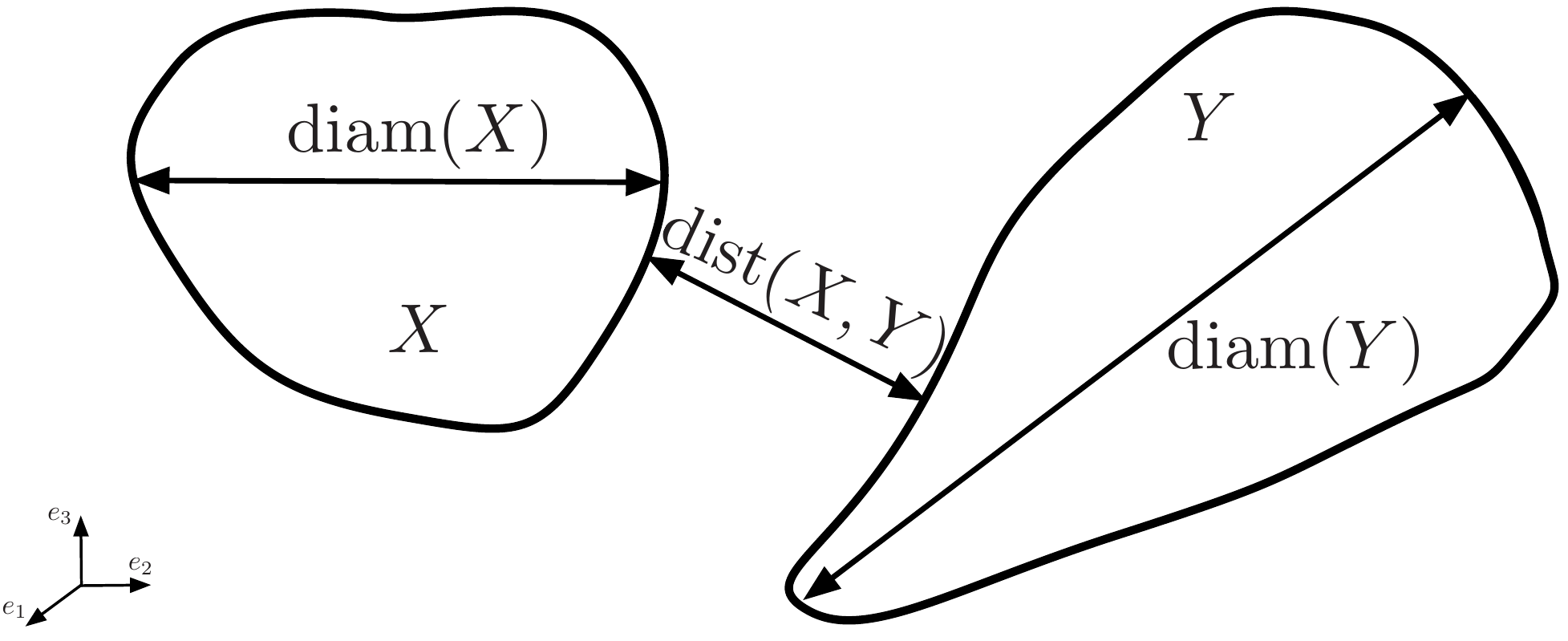} 
\end{tabular}
\end{center}
\caption{Definition of the distance and diameter.   }
\label{distances2}
\end{figure}

In other words, the exponential convergence of the Taylor expansion $s^r$ over $X$ and $Y$, with center $\x_0 \in X$ is constrained by a condition on $\x_0$, $X$ and $Y$. One can derive a sufficient condition, independent of $\x_0$, be observing that
\[
\max_{\x \in X} ||\x - \x_0|| \le \operatorname{diam}X \quad \mbox{and} \quad \operatorname{dist}(\x_0,Y) \ge \operatorname{dist}(X,Y)
\]
where   $\operatorname{dist}$ denotes the euclidian distance, i.e.
\begin{equation}
\mbox{dist}(X,Y)=\inf \{ ||\x-\y||, \ \x \in X, \ \y \in Y\}
\label{dist}
\end{equation}
and $\mbox{diam}$ denotes the diameter of a domain (Fig.~\ref{distances2}), i.e.
\begin{equation}
\mbox{diam}(X)=\sup \{ ||\x_1-\x_2||, \ \x_1, \ \x_2 \in X\}.
\label{diam}
\end{equation}

Let us  introduce the notion of $\eta$-admissibility.
   $X$ and $Y$ are said to be $\eta$-\emph{admissible} if $\operatorname{diam}X \le \eta \operatorname{dist}(X,Y)$ where $\eta>0$ is a parameter of the method. In this configuration, it follows that
\begin{equation}
\label{admi_cond_taylor}
\max_{\x \in X} ||\x - \x_0|| \le \operatorname{diam}X \le \eta \operatorname{dist}(X,Y) \le \eta \operatorname{dist}(\x_0,Y)  \quad \mbox{and so} \quad \gamma_x \le \frac{\eta}{c_2}.
\end{equation}
For asymptotically smooth kernels, we thus conclude that the exponential convergence is ensured if $X$ and $Y$ are $\eta$-\emph{admissible} and if $\eta$ is chosen such that $\eta < c_2$. We observe that in~\cite{bebendorf2008hierarchical}, the condition is more restrictive on $\eta$ since it reads $\eta \sqrt{3} < c_2$.

The Helmholtz Green's function  $G_{\kappa}(\x,\y) =\frac{\exp (i \kappa ||\x-\y||)}{4 \pi ||\x-\y||}$ is    not asymptotically smooth~\cite{banjai2008hierarchical}. On the other hand, noting that it holds
$
 G_{\kappa}(\x,\y)  =\exp (i \kappa ||\x-\y||)G_0(\x,\y) $
 and since $G_0(\x,\y) $ is   asymptotically smooth~\cite{bebendorf2008hierarchical,hackbusch2015hierarchical}, 
the Helmholtz Green's function  satisfies the following estimate~\cite{banjai2008hierarchical}: $\exists c_1,c_2$ and $\sigma \in \mathbb{N}_0$ (singularity degree of the kernel $G_0(.,.)$)
 such that $\forall z \in \{x_{\alpha},y_{\alpha}\}$, $\forall n \in \mathbb{N}_0$, $\forall \x \neq \y$
 \begin{equation}
 \label{est_Helmholtz}
 |\partial_z^n G_{\kappa}(\x,\y)| \le n! c_1 (1+\kappa ||\x-\y||)^n (c_2 ||\x-\y||)^{-n -\sigma}. 
\end{equation}
From inequality~(\ref{est_Helmholtz}), it is clear that in the so-called \emph{low-frequency regime}, i.e. when \linebreak $\kappa \ \max_{\x,\y \in X \times Y} (\kappa ||\x-\y||)$ is small, the  Helmholtz Green's function  behaves similarly to an asymptotically smooth kernel because 
$\kappa ||\x-\y||$ is uniformly small for all $\x,\y \in X \times Y$. 
There exists a lot of works on the $\mathcal{H}$-matrix representation~\cite{bebendorf2000approximation,borm2005hybrid} of asymptotically smooth kernels such as $G_0$ but much less attention  has been devoted to the case of the
3D Helmholtz equation. We can mention the  work~\cite{stolper2004computing} 
where a collocation approach is considered for problems up to $20 \ 480$ DOFs and where encouraging numerical results are presented.
 
We denote $G^r_{\kappa}$ the  Taylor expansion with $r:=| \{{\boldsymbol \alpha}\in \mathbb{N}^3_0:  |{\boldsymbol \alpha}| \le m\}|$ terms of the Helmholtz Green's function $G_{\kappa}$, i.e. 
\begin{equation}
G^r_{\kappa}(\x,\y)=\sum_{{\boldsymbol \alpha}\in \mathbb{N}^3_0, |{\boldsymbol \alpha}| \le m }   (\x-\x_0)^{\boldsymbol \alpha} \frac{1}{{\boldsymbol \alpha}!} \partial_x^{\boldsymbol \alpha}G_{\kappa}(\x_0,\y) + R^r_{\kappa} \quad \x \in X, \ \ \y \in Y
\label{taylor_exp}
\end{equation}
where $\x_0\in X$ is the centre of the expansion.
According to the previous discussion, we know that, for the asymptotically smooth kernel $G_0$, the asymptotic smoothness allows to prove the exponential convergence of the Taylor series if $X$ and $Y$ are $\eta$-\emph{admissible}. On the other hand for the Helmholtz kernel $G_{\kappa}$, 
the estimate~(\ref{est_Helmholtz}) now leads to
\[
|R^r_{\kappa} (\x,\y)|\le C' \sum_{\ell=m}^{\infty} \Big((1 + \kappa ||\x_0 -\y||)  \frac{ ||\x-\x_0||}{c_2 ||\x_0-\y||}\Big)^{\ell}.
\]
Thus to assure convergence of the series, $X$ and $Y$ must now  be chosen in such a way that
\begin{equation}
\gamma_{\kappa,x}:= \displaystyle    \ (1+\kappa ||\x_0- \y||)    
\frac{\max_{\x \in X}||\x-\x_0||}{c_2 ||\x_0-\y||} <1.
\label{admi_cond_taylor2}
\end{equation}
 From now on, we say  that $X$ and $Y$ are $\eta_{\kappa}$-\emph{admissible} if (\ref{admi_cond_taylor2}) is satisfied.
Provided that  the condition (\ref{admi_cond_taylor2}) holds, the remainder is bounded
\begin{equation}
\label{error_est_remain}
\begin{array}{rl}
|R^r_{\kappa}(\x,\y)|\le &\displaystyle{C'  
 \frac{\gamma_{\kappa,x}^m}{1-\gamma_{\kappa,x}}}.
 \end{array}
\end{equation}
Remark that in the low-frequency regime, one has the estimate $\gamma_{\kappa,x} \sim   
\frac{\max_{\x \in X}||\x-\x_0||}{c_2 ||\x_0-\y||} 
   $ such that the  $\eta_{\kappa}$-admissibility condition   is similar to the case of the asymptotically smooth kernel $s=G_0$. Within this framework, for higher frequencies, $\gamma_{\kappa,x}$ depends \emph{linearly} on the wavenumber $\kappa$ and thus on the circular frequency $\omega$.
  
  \begin{remark}
We do not discuss in this work on the choice of the optimal $\x_0$. It is very likely that the notions of $\eta$- and $\eta_{\kappa}$-admissibility can be improved by determining the optimal $\x_0$.
  \end{remark}
It is worth noting that a Taylor expansion with centre $\y_0\in Y$ instead of $\x_0\in X$  can be  performed. In that case, a similar error estimate is obtained:
\begin{equation}
\label{error_est_remain2}
\begin{array}{rl}
|R^r_{\kappa}(\x,\y)|\le &\displaystyle{C'  
 \frac{\gamma_{\kappa,y}^m}{1-\gamma_{\kappa,y}}}, \quad \mbox{ with } \quad \gamma_{\kappa,y}:= \displaystyle    \  (1+\kappa ||\y_0- \x||)
\frac{\max_{\y \in Y}||\y-\y_0||}{c_2 ||\y_0-\x||}. 
 \end{array}
\end{equation}

 \subsection{Theoretical estimate of the rank of 3D elastodynamics single-layer and double-layer potentials \label{taylor_elasto}}

 The results on the convergence of  the Taylor expansion of the 3D Helmholtz Green's function are easily extended to 3D elastodynamics. From~(\ref{elasto_U}) it is clear that the 3D elastodynamic Green's tensor is a combination of derivatives of the 3D Helmholtz Green's tensor. In particular, rewriting~(\ref{elasto_U}) component by component, we have
 \[
 ({\bf U}_{\omega})_{\alpha \beta}(\x,\y)=\frac{1}{\kappa_s^2 \mu} \Big((\delta_{\tau \sigma} \delta_{ \alpha \beta}-\delta_{ \beta \tau}\delta_{\alpha \sigma}) \frac{\partial}{\partial x_{\tau}} \frac{\partial}{\partial y_\sigma}G_{\kappa_s}(\x,\y)+   \frac{\partial}{\partial x_{\alpha}} \frac{\partial}{\partial y_{\beta}}G_{\kappa_p}(\x,\y) \Big)
 \] 
 \[
 ({\bf T}_{\omega})_{\alpha \beta}(\x,\y)= C_{\alpha \gamma \tau \sigma} \frac{\partial}{\partial y_{\sigma}}  \Big({\bf U}_{\omega}(\x,\y)\Big)_{\beta \tau} n_{\gamma}(\y)
 \] 
 where the Einstein summation convention is used.
%
Since $\kappa_p \le \kappa_s$, it follows that each component of the 3D elastodynamics Green's tensor satisfies an inequality similar to~(\ref{est_Helmholtz}) with $\kappa=\kappa_s$, i.e.
 \begin{equation}
 \label{est_elasto}
 |\partial_z^n ({\bf P})_{\alpha \beta}(\x,\y)| \le n! c_1 (1+\kappa_s ||\x-\y||)^n (c_2 ||\x-\y||)^{-n -\sigma}
\end{equation} 
where  ${\bf P}={\bf U}_{\omega}$ or ${\bf P}={\bf T}_{\omega}$. As a result, provided that (\ref{error_est_remain}) or~(\ref{error_est_remain2}) holds with $\kappa=\kappa_s$ the Taylor expansion of the 3D elastodynamic Green's tensor converges exponentially with convergence rate $\gamma_{\kappa_s,x}$ or $\gamma_{\kappa_s,y}$. However again two regimes can be distinguished.
In the  \emph{low-frequency} regime, the Taylor expansion behaves similarly to the case of the asymptotically smooth kernel $G_0$.
 The exponential convergence is ensured if $X$ and $Y$ are $\eta$-\emph{admissible} and if $\eta$ is chosen such that $\eta < c_2$. 
For higher frequencies, $\gamma_{\kappa_s,x}$ and $\gamma_{\kappa_s,y}$  defined in (\ref{error_est_remain}) and~(\ref{error_est_remain2}) depend linearly on the circular frequency $\omega$.

The definition of the  {$\eta_{\kappa}$-admissibility} for oscillatory kernels could be modified in order to keep the separation rank constant while the frequency increases. It is the option followed for example  in directional approaches~\cite{messner2012fast,delamotte2016etude}.  Another approach to avoid the rank increase is the concept  of $\mathcal{H}^2$-matrices~\cite{banjai2008hierarchical}.  
  The   option we consider  is to work with the   $\eta$-admissibility condition defined for the asymptotically smooth kernels $G_0$.
 In other words, we choose to  keep $\eta$ constant for all frequencies and use the admissibility condition defined by: 
\begin{equation}
\operatorname{diam}X \le \eta \operatorname{dist}(X,Y) \quad \mbox{ with } \eta<c_2.
\label{admi_G}
\end{equation}
Our goal is to determine what is the actual behavior of the algorithm in that configuration. Following the study of the remainder of the Taylor expansion, it is clear that the method will not be optimal since the separation rank to achieve a given accuracy will increase with the frequency. Nevertheless,  we are interested in determining the limit until which the low-frequency approximation is viable for oscillatory kernels such as the 3D Helmholtz Green's function or the 3D elastodynamic Green's tensors.

The relevant question to address  is then to determine the law followed by the growth of the numerical rank if the $\eta$-admissibility condition~(\ref{admi_G}) is used in the context of oscillatory kernels. 
 Importantly,  { the framework of } Taylor expansions is only used to illustrate the degenerate kernel in an abstract way  {and to determine the optimal admissibility conditions in the low and high-frequency regimes within this framework}. In practice   the ACA  {which is a purely algebraic method} is used after discretization.  {As a result, we will use linear algebra arguments to study the growth of the numerical rank.}
 As a matter of fact, one expects a small numerical rank provided that the entries of the matrix are approximated by a short sum of products.
  From the Taylor expansion
 (\ref{taylor_exp}) this can be achieved   if the discretization is accurate enough. 
 Consider two clusters of points $X=(\x_i)_{i \in \tau}$ and $Y=(\y_j)_{j \in \sigma}$ and the corresponding subblock $\mathbb{A}_{\tau \times \sigma}$. Then using  (\ref{taylor_exp}) we can build an approximation of $\mathbb{A}_{\tau \times \sigma}$ of the form $\mathbb{A}_{\tau \times \sigma} \approx \mathbb{B}_{\tau \times \sigma} =\mathbb{U}_{\tau } \mathbb{V}^*_{ \sigma}$
 where $\mathbb{U}_{\tau }$ corresponds to the part that depends on $\x$ and  $\mathbb{V}_{ \sigma}$ to the part that depends on $\y$. 
 
  In the following,   $\kappa$ denotes either  the wavenumber in the context of the 3D Helmholtz equation, or also the S-wavenumber $\kappa_s$ in the context of 3D elastodynamics.

 To start,  let us have a look at two clusters of points $X$ and $Y$ which are $\eta$-\emph{admissible} under~(\ref{admi_G}).
At the wavenumber $\kappa$ two cases may follow, whether condition~(\ref{admi_cond_taylor2}) is also satisfied or not.

Let's consider first that the condition~(\ref{admi_cond_taylor2}) is not satisfied at the wavenumber $\kappa$, i.e. if the frequency is too large, with respect to the diameters of the clusters, to consider that this block is $\eta_{\kappa}$-\emph{admissible}. From the study of the Taylor expansion, it follows that the numerical rank $r(\kappa)$ of the corresponding block   in the matrix at wavenumber $\kappa$ will be large and the algorithm will not be efficient.

Let's consider now the case where   $X$ and $Y$  are  such that the condition~(\ref{admi_cond_taylor2}) is also satisfied at the wavenumber $\kappa$. In that case, the numerical rank $r(\kappa)$  will be small. 

These two cases illustrate the so-called \emph{low-frequency} regime (latter case), i.e. when $(1+\kappa ||\x_0- \y||) \sim 1$ and \emph{high-frequency} regime (former case), i.e. when $(1+\kappa ||\x_0- \y||) \sim \kappa ||\x_0- \y||$.

 {In this work however, we are not interested by only one frequency but rather by a frequency range.
The key point is to remark that in the $\eta_{\kappa}$-admissibility condition, the important quantity is $\kappa ||\x_0- \y||$, i.e. the product between the wavenumber and  the diameter of the blocks (we recall that the $\eta$-admissibility~(\ref{admi_G}) relates the distances between blocks and the diameters of the blocks).
 Let's now consider the following configuration: at the wavenumber $\kappa$ all the blocks of a matrix are such that the $\eta$-admissibility includes the  $\eta_{\kappa}$-admissibility
while  at the  wavenumber $2\kappa$, it is not the case. This configuration represents the transition regime after the \emph{low-frequency} regime.
From  (\ref{error_est_remain}) and~(\ref{error_est_remain2}) it follows that  if the frequency  doubles, and as a result the wavenumber doubles from $\kappa$ to $\kappa'=2\kappa$,  an $\eta_{\kappa}$-\emph{admissible} block ${\mathbb A}_{\tau \times \sigma}$   should be decomposed into a $2 \times 2$ block matrix at wavenumber $\kappa'$ to reflect the reduction by a factor $2$ of the leading diameter. Thus the subblocks
  ${\mathbb A}_{\tau_i \times \sigma_j}$ are shown by the Taylor expansion to be of the order of the small rank $r(\kappa)$ at the wavenumber $\kappa'$, i.e.
 \[
{\mathbb A}_{\tau \times \sigma}=\left[ \begin{array}{cc}
{\mathbb A}_{\tau_1 \times \sigma_1} &{\mathbb A}_{\tau_1  \times \sigma_2}\\
{\mathbb A}_{\tau_2 \times \sigma_1}  & {\mathbb A}_{\tau_2 \times \sigma_2}\\
\end{array}
\right] \approx
\left[ \begin{array}{cc}
{\mathbb B}_{\tau_1 \times \sigma_1} &{\mathbb B}_{\tau_1 \times \sigma_2}\\
{\mathbb B}_{\tau_2 \times \sigma_1}  & {\mathbb B}_{\tau_2 \times \sigma_2}\\
\end{array}
\right]
\]
where each ${\mathbb B}_{\tau_i \times \sigma_j}$ ($1 \le i,j \le 2$) is a low-rank approximation of rank $r(\kappa)$ that can be written as ${\mathbb B}_{\tau_i \times \sigma_j}={\mathbb U}_{\tau_i \times \sigma_j}  {\mathbb S}_{\tau_i \times \sigma_j} {\mathbb V}_{\tau_i \times \sigma_j}^*$
and each $ {\mathbb S}_{\tau_i \times \sigma_j}$   is a $r(\kappa) \times r(\kappa)$ matrix.  
What is more,  the Taylor expansion (\ref{taylor_exp}) says that ${\mathbb U}_{\tau_i \times \sigma_j}  $ is independent of $ \sigma_j$ such that ${\mathbb U}_{\tau_i \times \sigma_j}  ={\mathbb U}_{\tau_i } $:
\[
{\mathbb A}_{\tau \times \sigma}\approx \left[ \begin{array}{cc}
{\mathbb U}_{\tau_1} & {\mathbb O}  \\
{\mathbb O} & {\mathbb U}_{\tau_2}\\
\end{array}
\right]
\left[ \begin{array}{cc}
 {\mathbb S}_{\tau_1 \times \sigma_1}  {\mathbb V}^*_{\tau_1 \times \sigma_1} &  {\mathbb S}_{\tau_1 \times \sigma_2} {\mathbb V}^*_{\tau_1 \times \sigma_2}\\
{\mathbb S}_{\tau_2 \times \sigma_1}  {\mathbb V}^*_{\tau_2 \times \sigma_1}   & {\mathbb S}_{\tau_2 \times \sigma_2}  {\mathbb V}^*_{\tau_2 \times \sigma_2}\\
\end{array}
\right].
\]
Since one can perform
the Taylor expansion also with respect to  $\y$, ${\mathbb A}_{\tau \times \sigma}$ is in fact of the form
\[
{\mathbb A}_{\tau \times \sigma}\approx \left[ \begin{array}{cc}
{\mathbb U}_{\tau_1} & {\mathbb O}  \\
{\mathbb O} & {\mathbb U}_{\tau_2}\\
\end{array}
\right]
\left[ \begin{array}{cc}
{\mathbb S}_{\tau_1 \times \sigma_1}  & {\mathbb S}_{\tau_1 \times \sigma_2}\\
{\mathbb S}_{\tau_2 \times \sigma_1}   & {\mathbb S}_{\tau_2 \times \sigma_2} \\
\end{array}
\right]
\left[ \begin{array}{cc}
 {\mathbb V}^*_{\sigma_1} & {\mathbb O}  \\
{\mathbb O}  &  {\mathbb V}^*_{\sigma_2}\\
\end{array}
\right]
\]
As a result,
   ${\mathbb A}_{\tau \times \sigma}$  is  numerically at most  of rank $2r(\kappa)$ at the wavenumber $\kappa'$. 
}  The same reasoning can be extended to higher frequencies  such that the maximum rank of all the blocks of a given matrix is at most increasing linearly with the frequency. 

 So far, we have assumed that all the blocks of a matrix satisfy the same  conditions. However due to the hierarchical structure of $\mathcal{H}$-matrix representations (i.e. subblocks in the matrix at different levels in the partition $\mathcal{P}$ corresponding to clusters with different diameters), it is very likely that  the blocks of a matrix at a given wavenumber will not be all in the \emph{low-frequency }regime nor all in the  \emph{high-frequency} regime. In particular assuming the frequency is not \emph{too high} and  provided the hierarchical subdivision is deep enough, there will always be some blocks in the  \emph{low-frequency} regime, i.e. blocks for which the rank is low numerically since the $\eta_{\kappa}$-\emph{admissibility} condition is satisfied.  Then for $\eta$-\emph{admissible} blocks at higher levels, the previous discussion on the transition regime can be applied  {by induction} such that the numerical  rank of these blocks is expected to grow linearly through the various levels. {Note that this linear increase depends on the choice of the number of points at the leaf level. The crux is to initialize the induction process: for a given choice, one has a certain limit for the transition regime. If one changes the number of points at the leaf level, the limit changes too. As a result, one cannot continue the process to infinity. To a given choice of number of points per leaf, there is a limit frequency.}

 From all these arguments, we conclude that if we consider a fixed geometry at a fixed frequency and we increase the density of discretization points per wavelength, and as a consequence the total number of discretization points,
the maximum rank of the blocks can be either small or large depending on the value of the product between the frequency and the diameter of the blocks (fixed by the geometry). However   we expect that this rank will remain  fixed if the density of discretization points per wavelength  increases: as a matter of fact, the Taylor expansion with $r$ given, is better and better approximated, so the rank should converge to $r$. This conclusion will be confirmed by numerical evidences in Section~7. On the other hand, if we consider a fixed geometry with a fixed density of points per S-wavelength, the maximum rank of the blocks is expected to grow linearly with the frequency until the \emph{high-frequency regime} is achieved. This transition regime is in some sense a \emph{pre-asymptotic regime} and we demonstrate numerically in Section~8 its existence.
    In the \emph{high-frequency} regime, the $\mathcal{H}$-matrix representation is expected to  {be suboptimal with a  rank rapidly increasing.}
    
 {This paper presents   an original algebraic argument to intuit the existence of this pre-asymptotic regime. It would be interesting in the future to compare this pre-asymptotic regime with  other more involved hierarchical  matrix based approaches proposed for oscillatory kernels, e.g. directional $\mathcal{H}^2$-matrices and with the multipole expansions used in the Fast Multipole Method. Another interesting but difficult question is to determine theoretically  the limits of the pre-asymptotic regime.}

\section{Implementation issues \label{impl_issues}} 

\subsection{Efficient implementation of the $\eta$-admissibility condition\label{impl_issues1}}
 A key tool in hierarchical matrix based methods is the efficient determination of \emph{admissible} blocks in the matrix. The first remark from an implementation point of view comes from the possibility to perform the Taylor expansion of 3D elastodynamic Green's tensors either with respect to $\x \in X=(\x_i)_{i \in \tau}$ or with respect to $\y \in Y=(\y_j)_{j \in \sigma}$.
When a Taylor expansion with respect to $\x$ is performed, 
the clusters  $X$ and $Y$ are said to be $\eta$-\emph{admissible} if $\operatorname{diam}X \le \eta \operatorname{dist}(X,Y)$.
 On the other hand, when a Taylor expansion with respect to $\y$ is performed, 
  the clusters $X$ and $Y$ are said to be $\eta$-\emph{admissible} if $\operatorname{diam}Y \le \eta \operatorname{dist}(X,Y)$.
Thus the admissibility condition at the discrete level becomes
\begin{equation}
\mbox{The block } \tau \times \sigma \mbox{ is  \emph{admissible}} \quad \mbox{ if } \mbox{min}\Big(\mbox{diam}(X), \mbox{diam}(Y) \Big) <  \hat{\eta} \times \mbox{dist}(X,Y)
\label{admissibility_cond}
\end{equation}
where $\hat{\eta}$ is a numerical parameter of the method discussed in the next section but subjected to some constraints highlighted in Section~\ref{taylor_elasto}.

\begin{figure}
\begin{center}
\begin{tabular}{cc}
\includegraphics[scale=0.35]{./distances} &\includegraphics[scale=0.35]{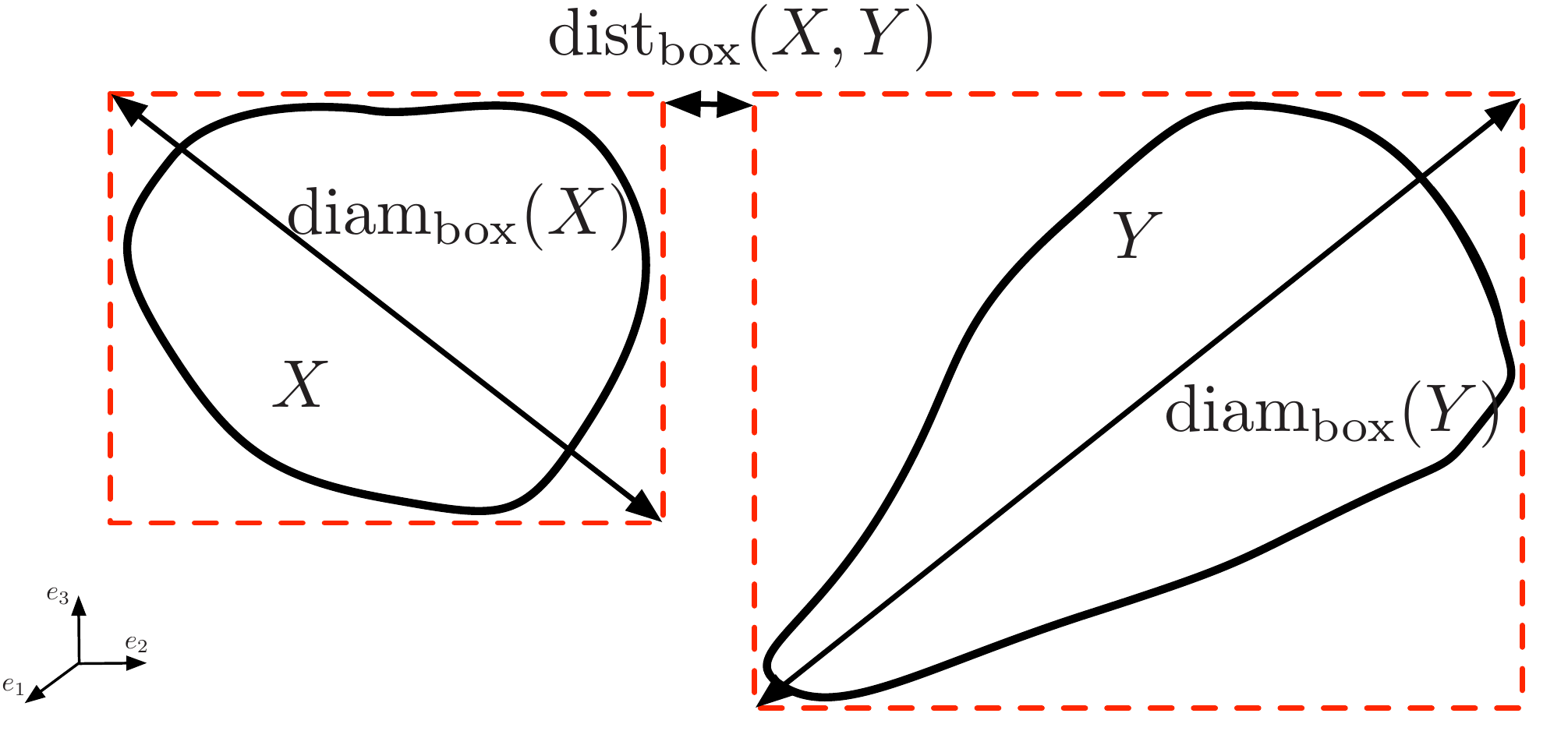} \\
(a) & (b)
\end{tabular}
\end{center}
\caption{Efficient implementation of the admissibility condition: (a) true condition and (b) more efficient approach using upper and lower bounds.   }
\label{distances}
\end{figure}

Since the computation of the diameter of a set $X$ with (\ref{diam}) is an expensive operation~\cite{hackbusch2015hierarchical}, it is replaced by computing the diameter of the bounding box with respect to the cartesian coordinates (Fig.~\ref{distances}b)
\[
\mbox{diam}(X) \le \Big( \sum_{\alpha=1}^3 (\max_{\x \in X} \x_{\alpha}- \min_{\x \in X} \x_{\alpha})^2 \Big)^{1/2}=: \mbox{diam}_{\operatorname{box}}(X)
\]
which is an upper bound. Similarly the distance between two sets is replaced by the distance between the closest faces of the bounding boxes
\[
\mbox{dist}(X,Y) \ge  \Big( \sum_{\alpha=1}^3 (\min_{\x \in X}\x_{\alpha} -\max_{\y \in Y}\y_{\alpha}  )^2+ (\min_{\y \in Y}\y_{\alpha} -\max_{\x \in X}\x_{\alpha}  )^2   \Big)^{1/2}=: \mbox{dist}_{\operatorname{box}}(X,Y)
\]
which is a lower  bound. 
The practical admissibility condition writes
\begin{equation}
\mbox{The block } \tau \times \sigma \mbox{ is  \emph{admissible}} \quad \mbox{ if } \mbox{min}\Big(\mbox{diam}_{\operatorname{box}}(X), \mbox{diam}_{\operatorname{box}}(Y) \Big) <  \hat{\eta} \times \mbox{dist}_{\operatorname{box}}(X,Y).
\label{admissibility_cond_bis}
\end{equation}

The admissibility conditions~(\ref{admissibility_cond})  or~(\ref{admissibility_cond_bis}) imply that 
 $X$ and $Y$ must be disjoint and that their
distance is related to their diameters. Also, one can check~\cite[Th. 6.16]{hackbusch2015hierarchical} that the sparsity pattern  
$C_{sp}$ is in fact bounded by a constant number that depends only on the geometry.

In addition, the choice of the parameter $\hat{\eta}$    influences    the shape of the \emph{admissible} blocks and  the number of \emph{admissible} blocks.
On the one hand, from the study of the Taylor expansion of the Green's function, it is clear that the smaller the parameter $\hat{\eta}$, the faster the remainder tends to zero.
On the other hand, a parameter $\hat{\eta}$ that tends to zero means that  only blocks of very small diameter (with respect to the distance) are \emph{admissible}.
The choice of a large $\hat{\eta}$ means in practice that  very elongated blocks are allowed. 
From an algorithmic point of view, our goal is to find a good compromise: i.e. not only  to have a large number of \emph{admissible}  blocks but also to ensure the convergence of the Taylor expansion. As a result, we want to choose $\hat{\eta}$ as large as possible.

 Since an upper bound of the diameter and a lower bound of the distance are used, it is clear that even though a theoretical bound on $\eta$ is found from the study of~(\ref{admi_G}),  larger values of $\hat{\eta}$ can be used and still lead to accurate  low rank block approximations in practice.
Then, the  distances and diameters are defined according to the cartesian coordinates cf.~(\ref{admissibility_cond_bis}). But depending on the choice of clustering method adopted, it is possible that the distance between some  blocks will be very small even though the blocks are disjoints. This remark weighs in favor of the use of a numerical parameter $\hat{\eta}$ larger than the one predicted by the theory.

%
%
%
%

\subsection{Choice of an acceptable parameter $\hat{\eta}$}

 To determine the adequate parameter $\hat{\eta}$, it is necessary to determine the constant $c_2$ in~(\ref{admi_cond_taylor2}) that gives the upper bound on $\eta$.
 This constant comes exclusively from the asymptotically smooth kernel $G_0$. It is easy to determine that
\[
 \forall z \in \{x_{\alpha},y_{\alpha}\} \quad \partial_z^1 G_0(\x,\y) \le \frac{1}{||\x-\y||^2} \quad  \mbox{ such that }  \quad 1\le c_1(c_2 ||\x-\y||)^{-2}.
\]
 Using Maple, it follows similarly that
\[
2\le c_1(c_2 ||\x-\y||)^{-3}, \ 4\le c_1(c_2 ||\x-\y||)^{-4},\ \frac{17}{2} \le c_1(c_2 ||\x-\y||)^{-5}     ,\ \frac{37}{2} \le c_1(c_2 ||\x-\y||)^{-6},
\]
\begin{align*}
  41 \le c_1(c_2 ||\x-\y||)^{-7}, \ 92 \le c_1(c_2 ||\x-\y||)^{-8}, \ \frac{5001}{24} \le c_1(c_2 ||\x-\y||)^{-9},\\
     \ \frac{3803}{8} \le c_1(c_2 ||\x-\y||)^{-10}, \ \frac{4363}{4} \le c_1(c_2 ||\x-\y||)^{-11}.
 \end{align*}
 In Table~\ref{reg_const}, we report the values of $c_1$ and $c_2$ obtained if a power regression up to the order $m$ is performed with $2 \le m \le 10$. We observe that $c_2$ is smaller than $0.5$ and decreases if $m$ is increased.   
 {As a result, in theory we should use  $\eta<c_2<0.5$.}
  Note that we have limited our study up to $m=10$ because $m\sim r^{1/3}$, where we recall that $r$ is the number of terms in the Taylor expansion~(\ref{taylor_exp}). We expect in practice to consider cases with a maximum numerical rank lower than $1000$ such that the current study is representative of the encountered configurations.
 \begin{table}[!htb]
 \begin{center}
 \begin{tabular}{c|ccccccccc}
  \hline
m & 2 & 3 &4 &5 & 6&7&8&9 &10\\
 \hline
$ c_1$ & $0.250$& $0.250$& $0.238$ & $0.226$& $0.216$& $0.206$& $0.198$& $0.191$& $0.184$\\
$ c_2 $&$0.500$&$0.500$&$0.490$&$0.483$&$0.476$&$0.470$&$0.465$&$0.461$&$0.458$ \\
 \hline
 \end{tabular}
 \end{center}
 \caption{Values of $c_1$ and $c_2$ obtained after a power regression for  values of $m$ between $2$ and $10$.}
  \label{reg_const}
  \end{table}

 {In Section~\ref{impl_issues1}, we have explained why  the use of   an upper bound of the diameter and a lower bound of the distance in the admissibility condition permits to  expect  in practice  good compression rates even though a larger value $\hat{\eta} >\eta$  is used.}
Indeed in~\cite{lize}, good numerical results for various geometries are obtained with  $\hat{\eta}=3$. This is the parameter we use in all our numerical examples for 3D elastodynamics.

\section{Behavior for low frequency elastodynamics  \label{fixed_freq}}

\subsection{Definition of the test problem  }
In Sections~\ref{fixed_freq} and~\ref{fixed_dens}, we consider the diffraction of   vertical incident plane P waves by a spherical cavity. The material properties are fixed to  $\mu=\rho=1$ and $\nu=1/3$.
We recall that in all the numerical examples, the binary tree $\mathcal{T}_I$ is built with a stopping criteria $N_{\mbox{leaf}}=100$ and the constant in the admissibility condition~(\ref{admissibility_cond}) is set to $\hat{\eta}=3$. Unless otherwise indicated, the required  accuracy is fixed to $\varepsilon_{\operatorname{ACA}}=10^{-4}$.
 In the following, we denote by $N$ the number of degrees of freedom such that $N=3N_c$.

\begin{figure}[!htb]
 \begin{center}
 \begin{tabular}{cc}
 \includegraphics[width=5cm]{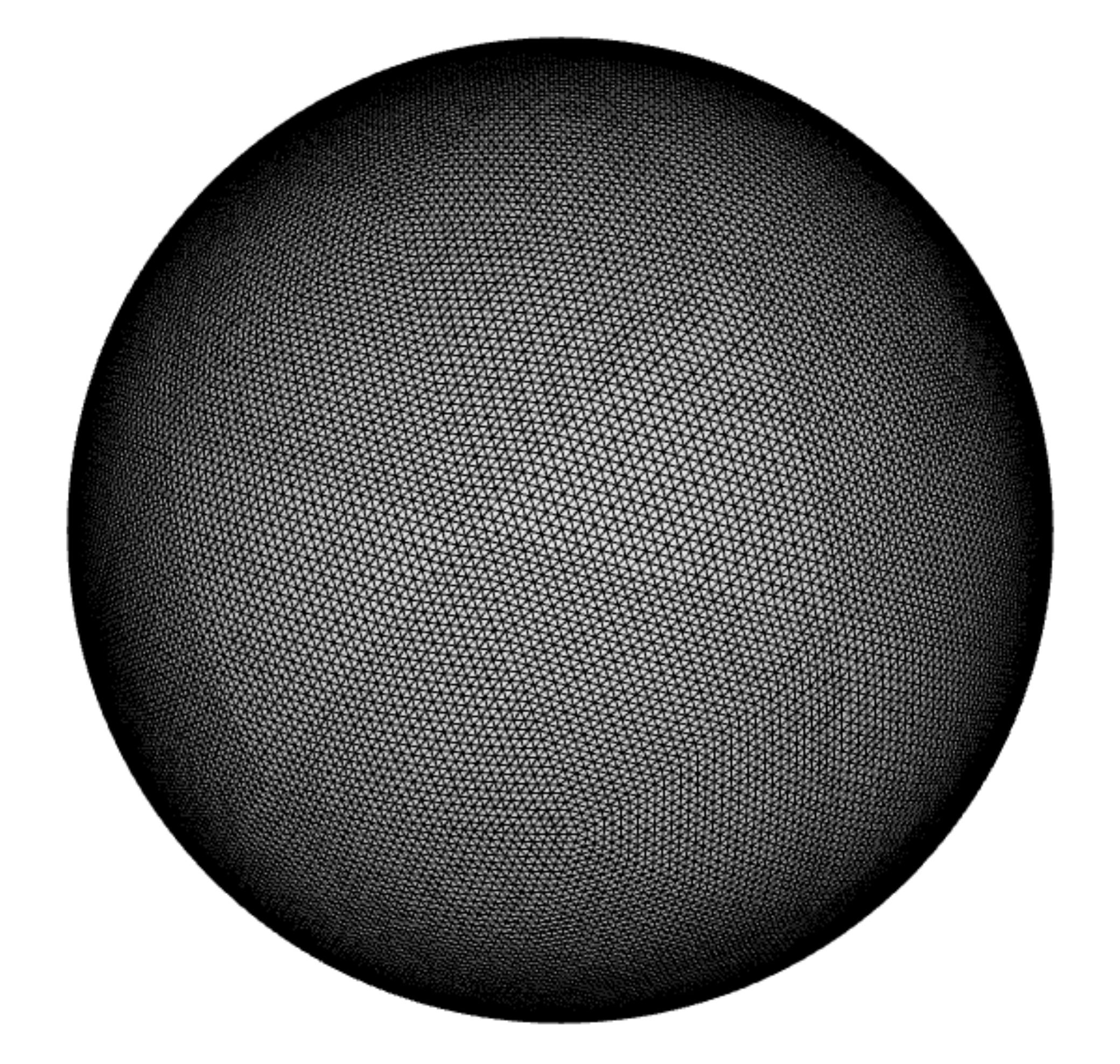}  
& \includegraphics[scale=0.45]{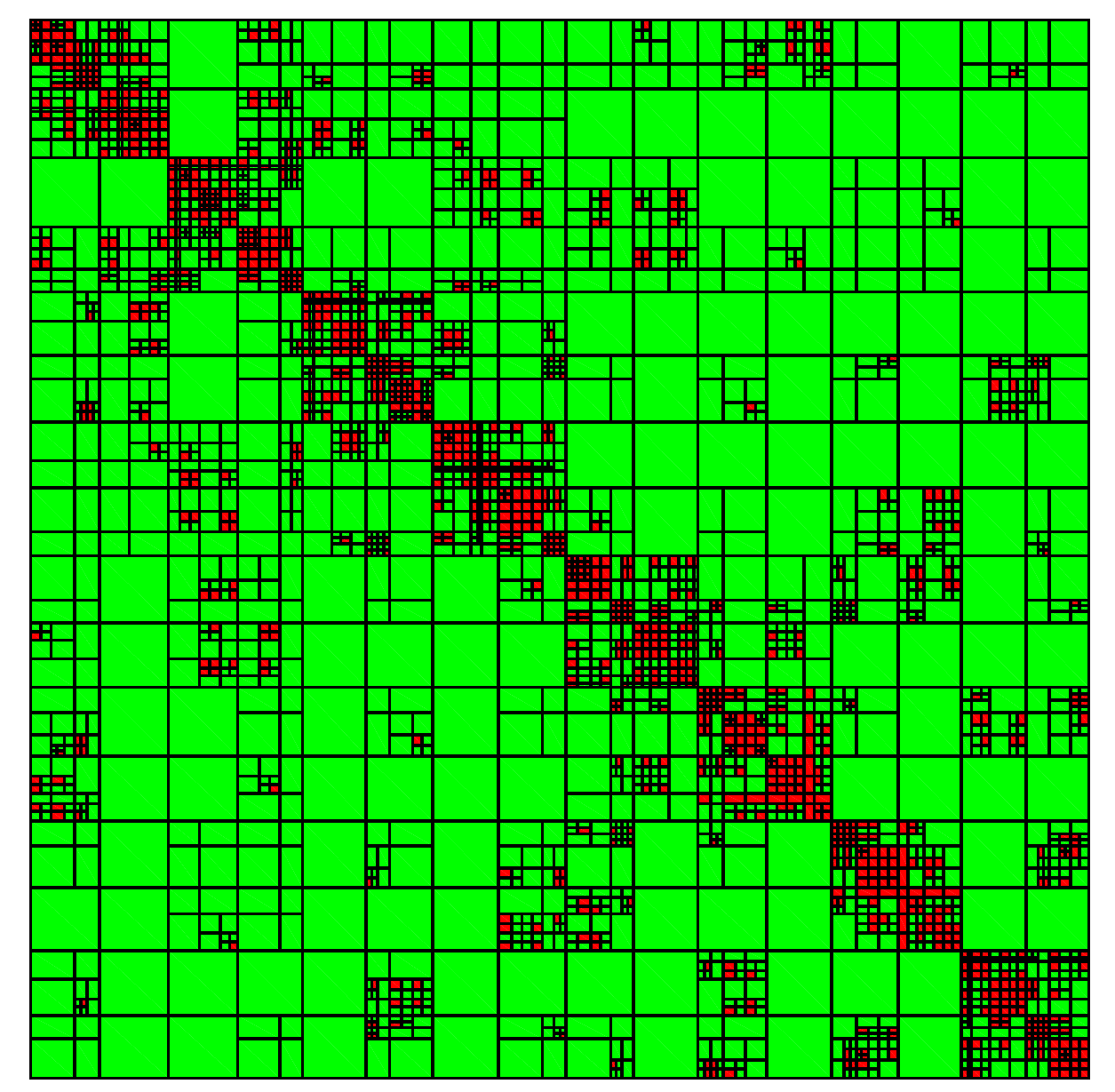} \\
 $N = 30 \ 726$ 
 \end{tabular}
 \end{center}
 \caption{Illustration of the expected repartition of the full and low-rank blocks for the sphere with $N= 30 \ 726$.} 
 \label{rep_mat}
 \end{figure}

The admissibility condition~(\ref{admissibility_cond}) depends only on the geometry of the domain. In Figure~\ref{rep_mat}, we illustrate the repartition of the expected low-rank (green  blocks) and full blocks (red blocks) for    the sphere geometry with $N= 30 \ 726$. 
 
 The goal of  Section~7 is to determine the behavior of the algorithm for a fixed frequency while the number of discretization points is increased (i.e. in the \emph{low-frequency} regime).
In the remainder of this section,  the circular frequency $\omega$ is  thus fixed.  As a result, the  density of points per S-wavelength increases as the number of discretization points increases. 

\subsection{Storage requirements for the single and double-layer potentials for a fixed frequency} 
In Section~\ref{theory}, we have seen that the memory requirements are expected to be of the order of \linebreak  $O(\max(r^{\max}_{\operatorname{ACA}}, N_{\operatorname{leaf}} )N \log_2 N)$. In addition, for a fixed frequency, $r^{\max}_{\operatorname{ACA}}$ is expected to be constant. 
In Tables~\ref{max_rank_frequency} and~\ref{max_rank_frequency2}, we report the maximum numerical rank observed among all the \emph{admissible} blocks for the case of a sphere for a fixed circular frequency ($\omega=3$ or $\omega=14$)
 and for various mesh sizes, for the single-layer and double-layer operators. The first rank corresponds to the numerical rank obtained by the partially pivoted ACA. The second rank corresponds to the numerical rank obtained after the recompression step (presented in Section~\ref{iter_solver}). The two ranks are seen to be almost constant while the density of points is increasing, as intuited by the study of the Taylor expansion. The numerical rank without the recompression step is much larger than after the optimisation step. The explanation  is that the partially pivoted ACA is a heuristic method, the rank is thus not the optimal one. The recompression step permits to recover the intrinsic numerical rank of the blocks. In addition, as expected the maximum numerical rank increases with the frequency. We will study in more details this dependence on the frequency in Section~8.
 Finally also as expected, since the two Green's tensors satisfy the inequality~(\ref{est_elasto}),   similar numerical ranks are observed for the single and double layer potentials. 

\begin{table}[!htb]
\begin{center}
\begin{tabular}{c|ccccccc}
\hline
$N$ & $7 \ 686$ & $30 \ 726$ & $122 \ 886$ & $183 \ 099$ & $ 490  \ 629$ & $763 \ 638$ &$ 985 \ 818$\\
\hline
$\omega=3$ & 63/39 & 72/39 & 75/39 & 69/37 & 75/39 & 69/40 & 78/39\\
$\omega=14$ & 99/73&105/75&108/76 &99/67& 114/76 & 111/76 & 111/76\\
\hline
\end{tabular}
\end{center}
\caption{Maximum numerical rank observed (before/after the  recompression step) for a fixed frequency, i.e. while increasing the density of points per S-wavelength, for $\varepsilon_{\operatorname{ACA}}=10^{-4}$ and the single-layer potential~(\ref{EFIE}).}
\label{max_rank_frequency}
\end{table}

\begin{table}[!htb]
\begin{center}
\begin{tabular}{c|ccccccc}
\hline
$N$ & $7 \ 686$ & $30 \ 726$ & $122 \ 886$ & $183 \ 099$ & $ 490  \ 629$ & $763 \ 638$ &$ 985 \ 818$\\
\hline
$\omega=3$ & 66/37 & 66/37 & 72/38 & 66/36  &  69/38 &   72/39& 75/39\\
$\omega=14$ & 102/70&117/74&117/74 & 102/66& 114/75 & 114/75&114/76\\
\hline
\end{tabular}
\end{center}
\caption{Maximum numerical rank observed (before/after the  recompression step) for a fixed frequency, i.e. while increasing the density of points per S-wavelength, for $\varepsilon_{\operatorname{ACA}}=10^{-4}$ and the double-layer potential~(\ref{MFIE}).}
\label{max_rank_frequency2}
\end{table}

\begin{figure}[!htbp]
 \begin{tabular}{cc}
\hspace*{-1cm} \includegraphics[scale=0.28]{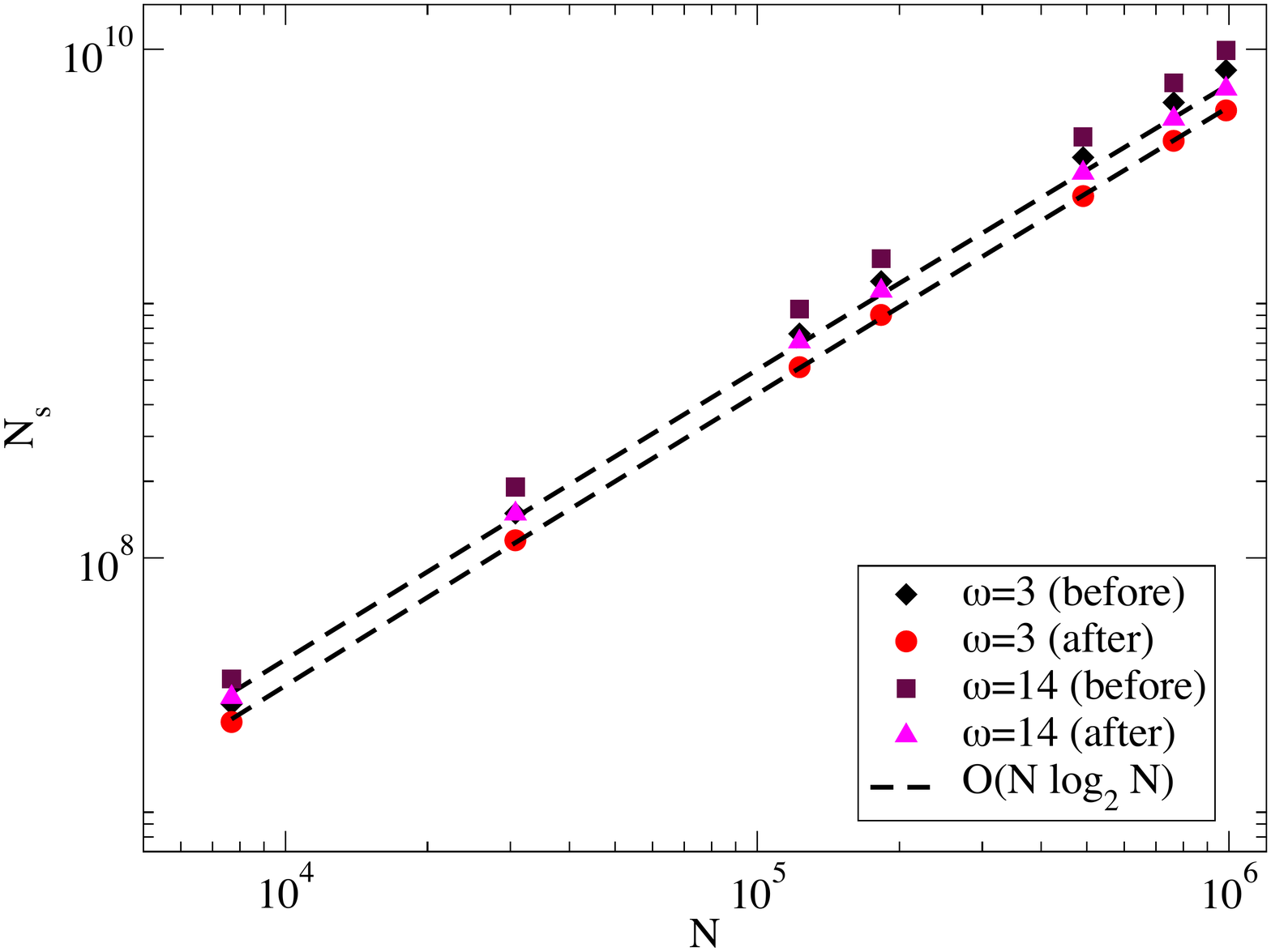} &  \includegraphics[scale=0.28]{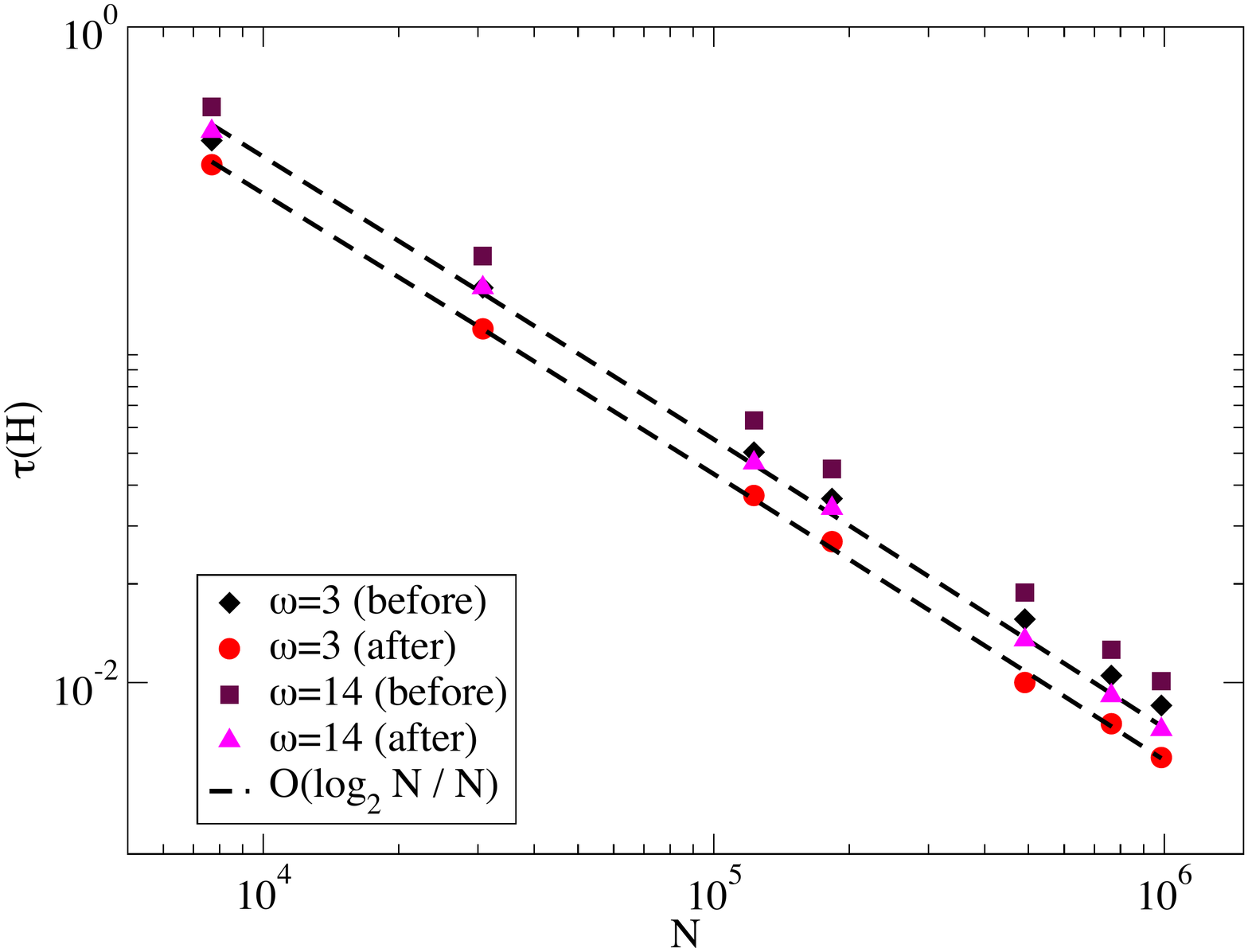} \\
({\bf a}) & ({\bf b})
 \end{tabular}
 \caption{ Observed and estimated   ({\bf a}) memory requirements $N_s$  and ({\bf b})  compression rate $\tau(\mathcal{H})$   with respect to the number of degrees of freedom $N$ for a fixed frequency ($\omega=3$ or $\omega=14$) for the single layer-operator and $\varepsilon_{\operatorname{ACA}}=10^{-4}$.}
 \label{graph_mem_frequency}
 \end{figure}

Ultimately we are not interested in the maximum numerical rank but rather on the memory requirements of the algorithm. In Figures~\ref{graph_mem_frequency}{\bf a-b}, 
  we report the memory requirements $N_s$ and the compression rate $\tau(\mathcal{H}):= N_s / N^{2}$ with respect to the number of degrees of freedom $N$ for $\varepsilon_{\operatorname{ACA}}=10^{-4}$ and the single-layer operator. Since the rank is constant, 
 $N_s$ is expected to be of the order of $N \log_2 N $ and   $\tau(\mathcal{H})$  of the order of $\log_2 N / N$. We observe for the two frequencies a very good agreement between the expected and observed complexities of the compression rate and memory requirements.
  Note that the storage is reduced by more than $95 \%$ as soon as $N \ge 10^{5}$.

 \subsection{Certification of the results of the direct solver} 
 
 While the previous subsection was dedicated to the study of the numerical efficiency of the method with respect to memory requirements, this subsection is more focused on the accuracy of the direct solver presented in Section~\ref{lu_solve}.
 Since the $\mathcal{H}$-LU factorization is based on various levels  of approximations, it is important to check the accuracy of the final solution. 
  In Section~\ref{esti_direct},
 we have derived a simple estimator to certify the results of the $\mathcal{H}$-LU direct solver. 
 
 Before we study  the accuracy, we recall that there are two important parameters in the method that  correspond to two sources of error: $\varepsilon_{\operatorname{ACA}}$ and $\varepsilon_{\operatorname{LU}}$.
 $\varepsilon_{\operatorname{ACA}}$ is tuned to the desired accuracy of the low-rank approximations of \emph{admissible} blocks, performed with the ACA described in Section~\ref{ACA_vect}. $\varepsilon_{\operatorname{LU}}$ is the parameter used  to adapt the   accuracy when recompressions are performed to modify  the format of a block   (Section~\ref{lu_solve}).  It is clear that the best accuracy that can be achieved by the direct solver is driven by  $\varepsilon_{\operatorname{ACA}}$. For this reason, the two parameters are set to the same value in the following numerical results.

 In Table~\ref{error_freq_single}, we report the value of the estimator  
  $\mathcal{I}(\delta_{\mathcal{H},F},\delta)$ and the observed normalized residual  $\frac{||{\bf b}-\mathbb{A}{\bf x}_0||_2}{ ||{\bf b} ||_2}$  if~(\ref{EFIE}) is solved with the $\mathcal{H}$-LU based direct solver.
 Two frequencies ($\omega=3$ and $\omega=14$) and four meshes are considered consistently with the previous subsection.
  When no compression is possible ($\tau(\mathcal{H})=1$ for $N=1 \ 926$) the estimated and observed residuals are, as expected, of the order of the machine accuracy. 
  For larger problems, the difference between the two residuals is not negligible. But importantly, when the required accuracy is reduced by two orders of magnitude, the two residuals follow the same decrease. In addition,  for a fixed value of  $\varepsilon_{\operatorname{ACA}}$, the accuracy of the solver and of the estimator does not deteriorate much while the number of degrees of freedom increases.
  
  To understand the lack of accuracy of our estimator, we add in Table~\ref{error_freq_single} the three components of the estimator: $\delta_{\mathcal{H},F}$ which corresponds to the level of error introduced by the low-rank approximation in the hierarchical representation of the system matrix (computed with the Frobenius norm), the ratio  
   $\frac{||{\bf x}_0||_2}{||{\bf b}||_2}$  and the normalized ratio $\frac{\delta}{||{\bf b}||_2}$ that accounts for the stability of the $\mathcal{H}$-LU factorization. On the one hand the contribution from  $\frac{\delta}{||{\bf b}||_2}$ is lower that the normalized residual that we try to estimate. On the other hand,   $ \delta_{\mathcal{H},F}$  is of the order of the required accuracy  $\varepsilon_{\operatorname{ACA}}$. But since this term is multiplied by the ratio     $\frac{||{\bf x}_0||_2}{||{\bf b}||_2}$ the part in the estimator accounting for the error  introduced by the low-rank approximations is overestimated.

  \begin{table}[!htb]
  \begin{center}
 \begin{tabular}{ccc|cc|ccc}
 \hline
 $\#$ DOFs & $\omega$ & $\varepsilon_{\operatorname{ACA}}=\varepsilon_{\operatorname{LU}}$ & $\mathcal{I}(\delta_{\mathcal{H},F},\delta)$ & $\frac{||{\bf b}-\mathbb{A}{\bf x}_0||_2}{ ||{\bf b} ||_2}$ & $\delta_{\mathcal{H},F}$ & $\frac{||{\bf x}_0||_2}{||{\bf b}||_2}$ & $\frac{\delta}{||{\bf b}||_2}$ \\
 \hline
 $1 \ 926$ & $3$ & $10^{-4}$ & $3.25 \ 10^{-15}$& $3.64 \ 10^{-15}$ & $9.16 \ 10^{-17}$ & $6.40$ & $2.66 \ 10^{-15}$ \\
 $7 \ 686$ & $3$ & $10^{-4}$ &$4.64 \ 10^{-4}$ & $9.94 \ 10^{-6}$ & $7.21 \ 10^{-5}$ & $6.37$ & $4.62 \ 10^{-6}$\\
 $30 \ 726$ & $3$ & $10^{-4}$ & $6.37 \ 10^{-4}$ & $1.18 \ 10^{-5}$ & $9.93 \ 10^{-5}$ & $6.36$ & $5.53 \ 10^{-6}$\\
  $122\ 886$ & $3$ & $10^{-4}$ &$8.14 \ 10^{-4}$ & $1.47 \ 10^{-5}$ & $1.27 \ 10^{-4}$ & $6.36$ & $5.21 \ 10^{-6}$\\ 
 $1 \ 926$ & $3$ & $10^{-6}$ & $3.25 \ 10^{-15}$& $3.64 \ 10^{-15}$ & $9.16 \ 10^{-17}$ & $6.40$ & $2.66 \ 10^{-15}$ \\
 $7 \ 686$ & $3$ & $10^{-6}$ &$5.07  \ 10^{-6}$ & $  7.35 \ 10^{-8}$ & $  7.92 \ 10^{-7}$ & $6.37$ & $ 2.48 \ 10^{-8}$\\
 $30 \ 726$ & $3$ & $10^{-6}$ & $ 7.44 \ 10^{-6}$ & $  9.38 \ 10^{-8}$ & $  1.16\ 10^{-6}$ & $6.36 $ & $ 2.93 \ 10^{-8}$\\  
   $122\ 886$ & $3$ & $10^{-6}$ & $9.25 \ 10^{-6}$ & $9.82 \ 10^{-8}$ & $1.45 \ 10^{-6}$& $6.36 $ & $2.76 \ 10^{-8}$ \\
 \hline
  $1 \ 926$ & $14$ & $10^{-4}$ & $5.69 \ 10^{-15}$& $4.69 \ 10^{-15}$ & $7.16 \ 10^{-17}$ & $29.29$ & $3.59 \ 10^{-15}$ \\
   $7 \ 686$ & $14$ & $10^{-4}$ &$3.89 \ 10^{-3}$ & $1.04 \ 10^{-4}$ & $1.36 \ 10^{-4}$ & $28.15$ & $4.73 \ 10^{-5}$\\
    $30 \ 726$ & $14$ & $10^{-4}$ & $5.49 \ 10^{-3}$ & $1.29 \ 10^{-4}$ & $1.95 \ 10^{-4}$ & $27.88$ &  $5.48 \ 10^{-5}$\\
      $122\ 886$ & $14$ & $10^{-4}$ & $6.24 \ 10^{-3}$& $1.49 \ 10^{-4}$&$2.22 \ 10^{-4}$& $27.82$ & $5.86 \ 10^{-5}$\\
     $1 \ 926$ & $14$ & $10^{-6}$ & $5.69 \ 10^{-15}$& $4.69 \ 10^{-15}$ & $7.16 \ 10^{-17}$ & $29.29$ & $3.59 \ 10^{-15}$ \\
     $7 \ 686$ & $14$ & $10^{-6}$ & $4.59 \ 10^{-5}$ & $9.15 \ 10^{-7}$ & $1.62 \ 10^{-6}$ & $28.15$ & $2.03 \ 10^{-7}$ \\  
       $30 \ 726$ & $14$ & $10^{-6}$ & $5.95 \ 10^{-5}$ & $1.15 \ 10^{-6}$ & $2.12 \ 10^{-6}$ &  $27.88$ &  $2.44 \ 10^{-7}$\\  
           $122\ 886$ & $14$ & $10^{-6}$ &$6.42 \ 10^{-5}$ & $1.26 \ 10^{-6}$& $2.30 \ 10^{-6}$& $27.82$ & $2.52 \ 10^{-7}$\\
 \hline
 \end{tabular}
 \caption{Accuracy of the direct solver for a fixed frequency ($\omega=3$ or $\omega=14$)  for the single-layer potential.}
 \label{error_freq_single}
 \end{center}
 \end{table}


 To accelerate the computation of the estimator, we have used the Frobenius norm. To check the effect of the use of this norm, we report in Table~\ref{error_freq_single2}
the values of the estimator $\mathcal{I}(\delta_{\mathcal{H},2},\delta)$ i.e. when the 2-norm is used. We consider only  moderate size problems due to the need of the computation of the singular value decomposition of the complete matrix. As expected, the estimator  $\mathcal{I}(\delta_{\mathcal{H},2},\delta)$ is much more accurate than  $\mathcal{I}(\delta_{\mathcal{H},F},\delta)$ but its cost is prohibitive. A good compromise consists in checking only the two main components of the estimator: $ \delta_{\mathcal{H},F}$ and $\frac{\delta}{||{\bf b}||_2}$ to check that the level of error  introduced in the different parts of the solver is consistent with the user parameter $\varepsilon_{\operatorname{ACA}}=\varepsilon_{\operatorname{LU}}$.

 \begin{table}
  \begin{center}
 \begin{tabular}{ccc|ccc|cc}
 \hline
 $\#$ DOFs & $\omega$ & $\varepsilon_{\operatorname{ACA}}=\varepsilon_{\operatorname{LU}}$ & $\mathcal{I}(\delta_{\mathcal{H},F},\delta)$ & $\mathcal{I}(\delta_{\mathcal{H},2},\delta)$ & $\frac{||{\bf b}-\mathbb{A}{\bf x}_0||_2}{ ||{\bf b} ||_2}$ & $\delta_{\mathcal{H},F}$ &  $\delta_{\mathcal{H},2}$  \\
 \hline
 $7 \ 686$ & $3$ & $10^{-4}$ &$4.64 \ 10^{-4}$ & $7.91 \ 10^{-5}$ & $9.94 \ 10^{-6}$ & $7.21 \ 10^{-5}$  & $1.17 \ 10^{-5}$ \\
 $30 \ 726$ & $3$ & $10^{-4}$ & $6.37 \ 10^{-4}$ &$5.46 \ 10^{-5}$ & $1.18 \ 10^{-5}$ & $9.93 \ 10^{-5}$ & $7.72 \ 10^{-6}$ \\ 
 $7 \ 686$ & $3$ & $10^{-6}$ &$5.07  \ 10^{-6}$ & $7.32 \ 10^{-7}$& $  7.35 \ 10^{-8}$ & $  7.92 \ 10^{-7}$ & $1.11 \ 10^{-7}$ \\
 $30 \ 726$ & $3$ & $10^{-6}$ & $ 7.44 \ 10^{-6}$ &$8.50 \ 10^{-7}$ & $  9.38 \ 10^{-8}$ & $  1.16\ 10^{-6}$ &  $1.29 \ 10^{-7}$\\   
 \hline
   $7 \ 686$ & $14$ & $10^{-4}$ &$3.89 \ 10^{-3}$ &$6.32 \ 10^{-4}$& $1.04 \ 10^{-4}$ & $1.36 \ 10^{-4}$ &  $2.08 \ 10^{-5}$\\
    $30 \ 726$ & $14$ & $10^{-4}$ & $5.49 \ 10^{-3}$ &$9.50 \ 10^{-4}$ & $1.29 \ 10^{-4}$ & $1.95 \ 10^{-4}$ &  $3.21 \ 10^{-5}$\\
     $7 \ 686$ & $14$ & $10^{-6}$ & $4.59 \ 10^{-5}$ &$1.07 \ 10^{-5}$& $9.15 \ 10^{-7}$ & $1.62 \ 10^{-6}$ & $3.72 \ 10^{-7}$  \\  
       $30 \ 726$ & $14$ & $10^{-6}$ & $5.95 \ 10^{-5}$ &$1.20 \ 10^{-5}$ & $1.15 \ 10^{-6}$ & $2.12 \ 10^{-6}$ & $4.20 \ 10^{-7}$  \\  
 \hline
 \end{tabular}
 \caption{Accuracy of the direct solver for a fixed frequency ($\omega=3$ or $\omega=14$)  for the single-layer potential when the 2-norm is used instead of the Frobenius norm.}
 \label{error_freq_single2}
 \end{center}
 \end{table}

 In Table~\ref{error_freq_double}, we  consider now the case where the double-layer potential based on equation~(\ref{MFIE})  is used and solved with the $\mathcal{H}$-LU  direct solver.
  It is clear that the behaviors of the single and double layer potentials are similar also in terms of numerical accuracy.    
    \begin{table}[!htb]
  \begin{center}
 \begin{tabular}{ccc|cc|ccc}
 \hline
 $\#$ DOFs & $\omega$ & $\varepsilon_{\operatorname{ACA}}=\varepsilon_{\operatorname{LU}}$ & $\mathcal{I}(\delta_{\mathcal{H},F},\delta)$ & $\frac{||{\bf b}-\mathbb{A}{\bf x}_0||_2}{ ||{\bf b} ||_2}$ & $\delta_{\mathcal{H},F}$ & $\frac{||{\bf x}_0||_2}{||{\bf b}||_2}$ & $\frac{\delta}{||{\bf b}||_2}$ \\
 \hline
  $1 \ 926$ & $14$ & $10^{-4}$ & $3.35 \ 10^{-15}$ & $3.28 \ 10^{-15}$ & $5.22 \ 10^{-16}$ & $1.26$ & $2.69 \ 10^{-15}$  \\
   $7 \ 686$ & $14$ & $10^{-4}$ & $2.49 \ 10^{-3}$ & $7.16 \ 10^{-5}$ & $1.94 \ 10^{-3}$ & $1.27$ & $2.58 \ 10^{-5}$  \\
    $30 \ 726$ & $14$ & $10^{-4}$ &$3.18 \ 10^{-3}$ & $6.92 \ 10^{-5}$ & $2.43 \ 10^{-3}$ & $1.30$ & $3.05 \ 10^{-5}$  \\
      $122\ 886$ & $14$ & $10^{-4}$ & $9.23 \ 10^{-3}$ & $5.70 \ 10^{-3}$ & $2.66 \ 10^{-3}$ & $1.33$ & $5.70 \ 10^{-3}$ \\
     $1 \ 926$ & $14$ & $10^{-6}$ &   $3.35 \ 10^{-15}$ & $3.28 \ 10^{-15}$ & $5.22 \ 10^{-16}$ & $1.26$ & $2.69 \ 10^{-15}$  \\
     $7 \ 686$ & $14$ & $10^{-6}$ & $2.43 \ 10^{-5}$ & $4.77 \ 10^{-7}$ & $1.90 \ 10^{-5}$ & $1.27$ & $9.85 \ 10^{-8}$   \\  
       $30 \ 726$ & $14$ & $10^{-6}$ & $3.31 \ 10^{-5}$ & $5.25 \ 10^{-7}$ & $2.54 \ 10^{-5}$ & $1.30$ & $1.15 \ 10^{-7}$ \\  
           $122\ 886$ & $14$ & $10^{-6}$ & $5.74 \ 10^{-3}$ & $5.70 \ 10^{-3}$ & $3.17 \ 10^{-5}$ & $1.33$ & $5.70 \ 10^{-3}$ \\  
 \hline
 \end{tabular}
 \caption{Accuracy of the direct solver for a fixed frequency ($\omega=14$)  for the double-layer potential.}
 \label{error_freq_double}
 \end{center}
 \end{table}
 
  All these numerical experiments confirm the theoretical study of the Taylor expansion, which yields that $\mathcal{H}$-matrix based solvers are very efficient tools for the \emph{low-frequency} elastodynamics.

\section{Towards highly oscillatory elastodynamics   \label{fixed_dens}}
The goal of  Section~8 is to determine the behavior of the algorithm for a given fixed density of points per S-wavelength  while the number of discretization points is increased (i.e. in the \emph{higher frequency} regime).
The circular frequency $\omega$ is thus fitted to the mesh length to keep the density of points per S-wavelength fixed for the different meshes. In this case, an increase of the number of degrees of freedom $N$ leads to an increase of the size of the body to be approximated. In Table~\ref{table_corr} we report the   number of degrees of freedom, corresponding  circular frequencies 
and number of S-wavelengths spanned in the sphere diameter used in the examples. 
 
 \begin{table}[!htbp]
 \begin{center}
 \begin{tabular}{c|ccccccc}
 $\omega$ & 8.25 & 16.5 & 33 & 40 & 66 & 84 & 92\\
 \hline
 $N$& $7 \ 686$ & $30 \ 726$ & $122 \ 886$ & $183 \ 099$ & $ 490  \ 629$ & $763 \ 638$ &$ 985 \ 818$\\
$n_{ \lambda_s}$& 2.6 & 5.2& 10.5& 12.6& 21 &26.6& 29.1\\
\hline
 \end{tabular}
 \end{center}
 \caption{Number of degrees of freedom,  corresponding number of wavenumbers spanned in the sphere diameter $D$ ($D=n_{ \lambda_s} \lambda_s$) on the different meshes of spheres and corresponding circular frequencies used in the examples of Section~8 (i.e for a fixed density of points per S-wavelength).}
 \label{table_corr}
 \end{table}

\subsection{Storage requirements for the single and double-layer potentials for a fixed density of points per S-wavelength}

\begin{figure}[!htbp]
 \begin{tabular}{cc}
\hspace*{-1cm} \includegraphics[scale=0.28]{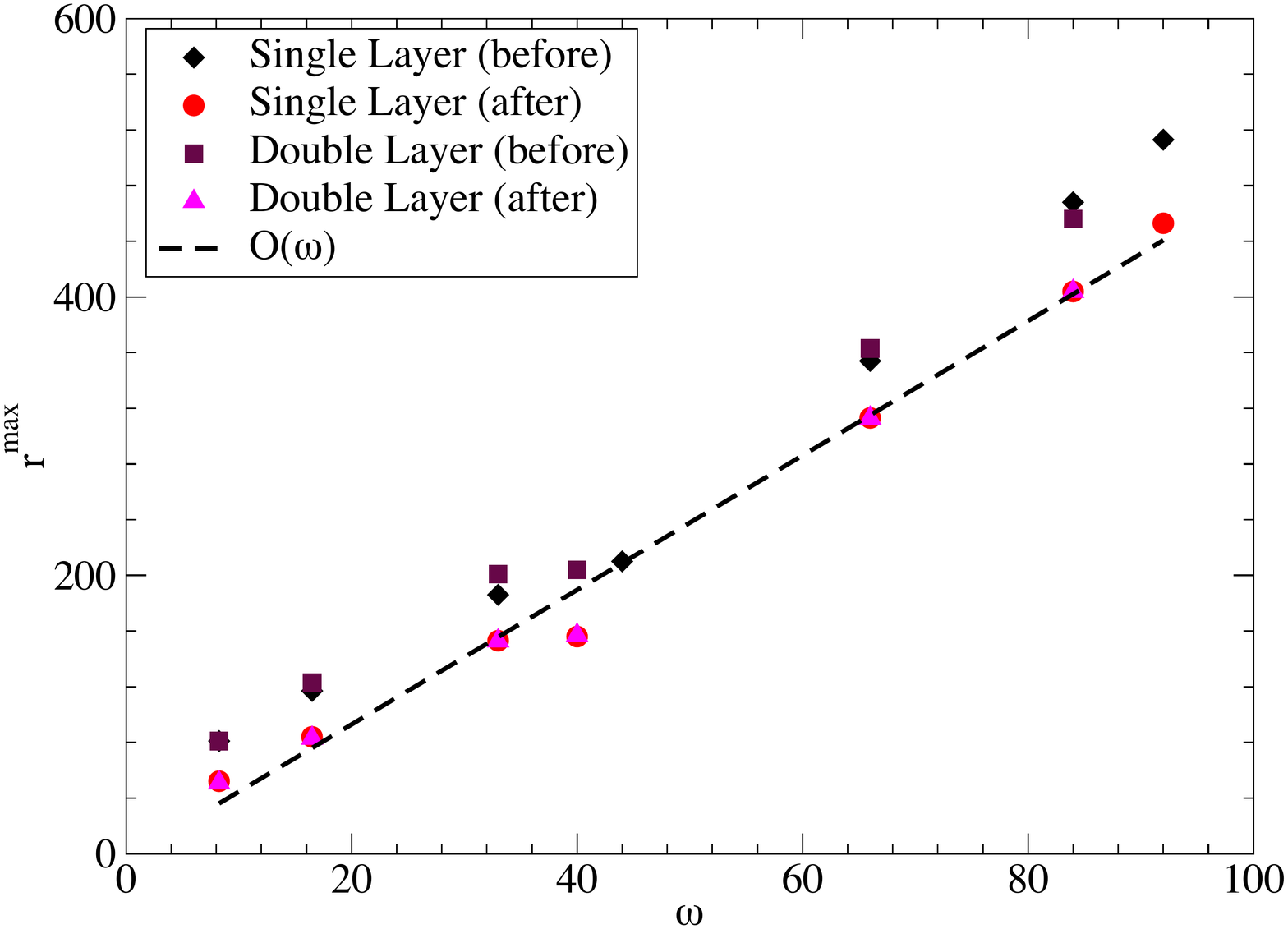} &  \includegraphics[scale=0.28]{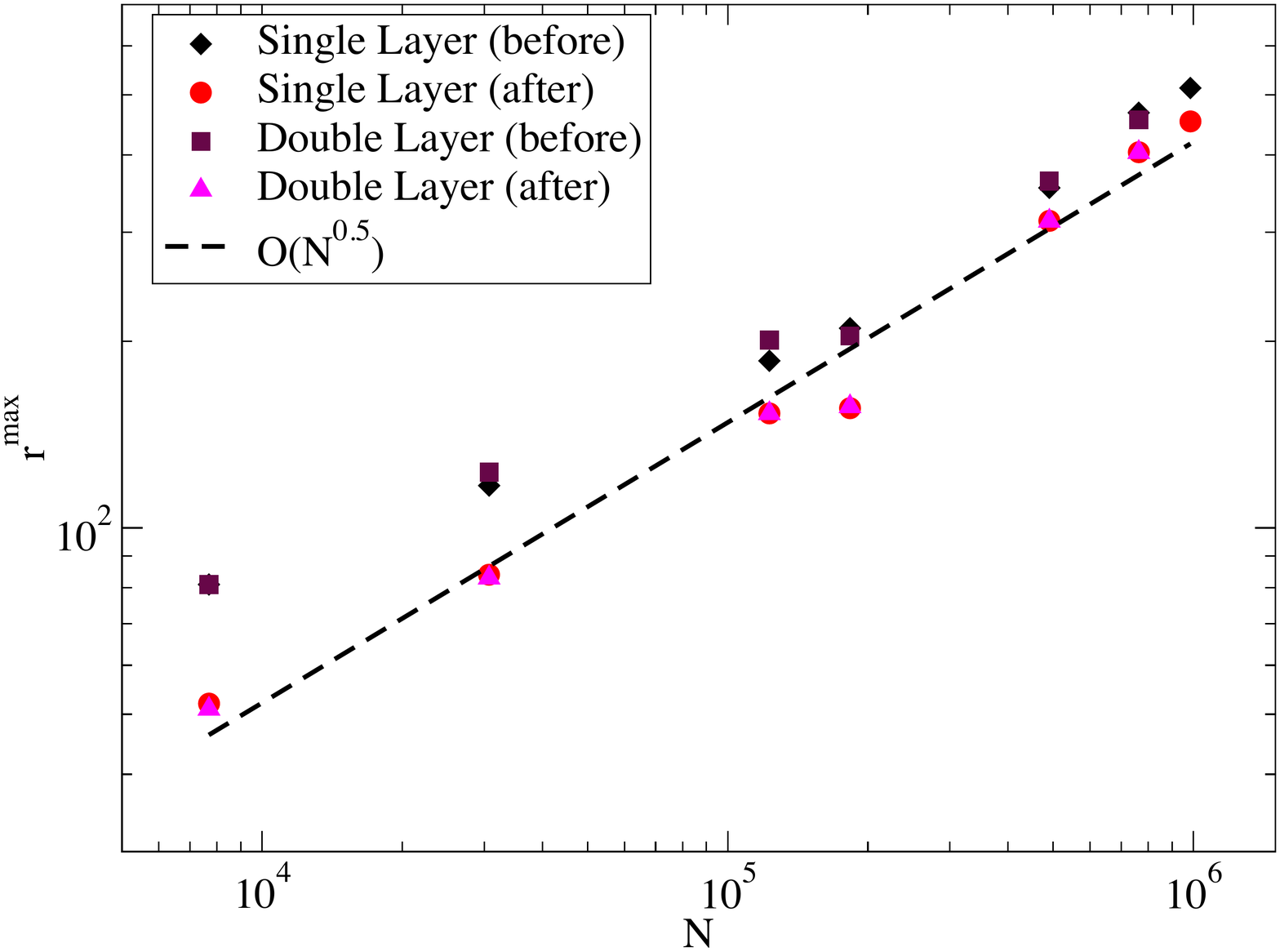} \\
{\bf (a)} & {\bf (b)}
 \end{tabular}
 \caption{Maximum numerical rank observed (before/after the  recompression step) for a fixed density of points per wavelength, i.e. while increasing the circular frequency with the number of discretization points, for $\varepsilon_{\operatorname{ACA}}=10^{-4}$: {\bf (a)}  with respect to the circular frequency and {\bf (b)}  with respect to the number of degrees of freedom.}
 \label{rank_fixed_freq}
 \end{figure}

In Section~\ref{taylor_elasto}, it has been shown that for a fixed geometry with a fixed density of points per S-wavelength, the maximum rank of the $\eta$-\emph{admissible} blocks is expected to grow linearly with the frequency until the \emph{high-frequency regime} is achieved. The goal of this section is to validate the existence of this pre-asymptotic regime.
In Figure~\ref{rank_fixed_freq}{\bf a}, we report the maximum numerical rank observed among all the $\hat{\eta}$-\emph{admissible} blocks before and after the recompression step with respect to the circular frequency $\omega$,  for both the single and double layer potentials. In accordance to the study of the Taylor expansion, the growth of the numerical rank after the recompression step (i.e. the intrinsic rank of the blocks) is seen to be linear with respect to the circular frequency. The behavior of the single layer and double layer operators are again very similar.
 
 In practice however, we are interested in the complexity with respect to the number of degrees of freedom. In this section, the number of discretization points is fitted to the circular frequency in order to keep 10 discretization points per S-wavelength. Since the density of points is fixed, the mesh size is  $h=O(\lambda_S)=O(\omega^{-1})$. In addition since the mesh for the BEM is  a surface mesh, the mesh size is $h=O(N^{-2})$. Hence,   the maximum numerical rank is expected to be of the order of $O(N^{1/2})$. In Figure~\ref{rank_fixed_freq}{\bf b}, we report the maximum numerical rank observed among all the $\hat{\eta}$-\emph{admissible} blocks before and after the recompression step with respect to the number of DOFs $N$,  for both the single and double layer potentials. The growth of the numerical rank after the recompression step (i.e. the intrinsic rank of the blocks) is indeed seen to be  of the order of $O(N^{1/2})$. 

\begin{figure}[!htbp]
\begin{tabular}{cc}
\includegraphics[scale=0.28]{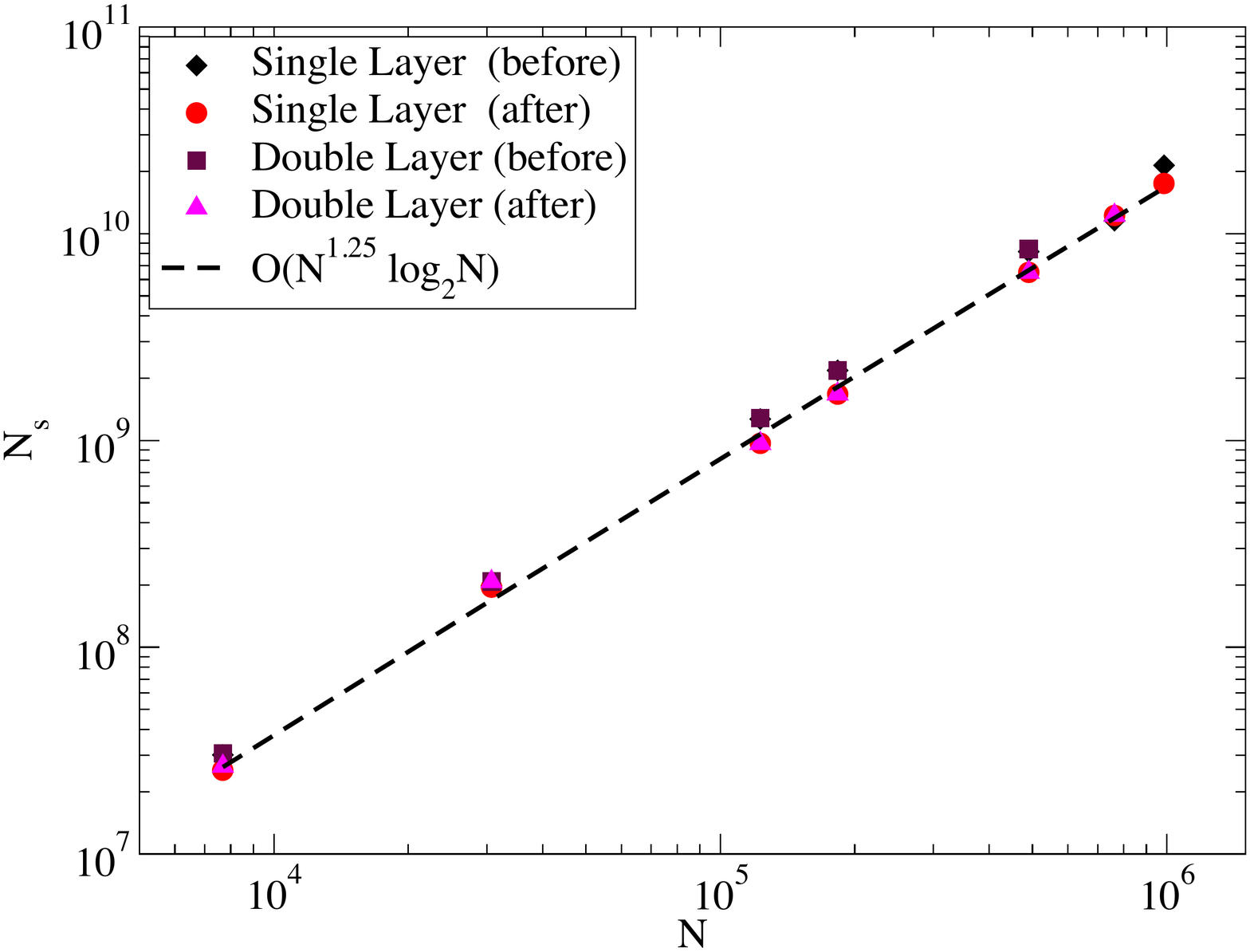}  & \includegraphics[scale=0.28]{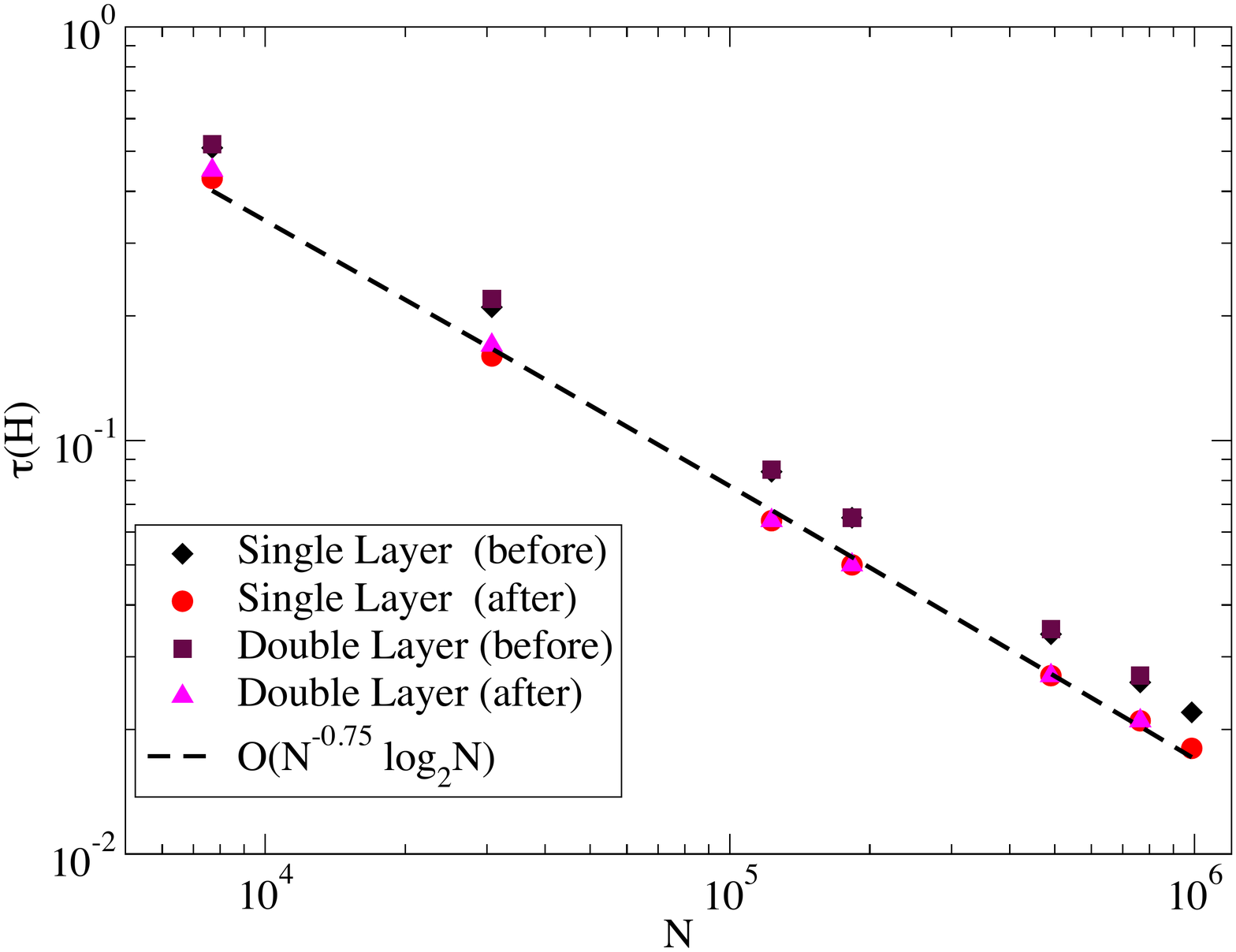} \\
{\bf (a)} & {\bf (b)}
 \end{tabular}
 \caption{Observed and estimated  ({\bf a}) memory requirements $N_s$     and ({\bf b})  compression rate $\tau(\mathcal{H})$ 
  with respect to the number of degrees of freedom $N$  
 for a fixed density of points per S-wavelength for the single and double layer-operators and $\varepsilon_{\operatorname{ACA}}=10^{-4}$.}
 \label{mem_fixed_freq}
 \end{figure}

 Ultimately we are not interested in the maximum numerical rank but rather on the memory requirements of the algorithm. In Figures~\ref{mem_fixed_freq}{\bf a-b}, 
  we report the memory requirements $N_s$ and compression rate $\tau(\mathcal{H})$,  with respect to the number of degrees of freedom $N$ for $\varepsilon_{\operatorname{ACA}}=10^{-4}$. Since the rank is of the order of $O(N^{1/2})$, 
$N_s$ is expected to be of the order of $N^{3/2} \log_2 N $ and  $\tau(\mathcal{H})$  of the order of $\log_2 N / N^{1/2}$. In practice,    observed complexities are lower than the expected one. The reason is that the estimation of the memory requirements gives only an upper bound based on the estimation of the maximum numerical rank. But, all the \emph{admissible} blocks do not have the same numerical rank such that the complexity is lower than the estimated one.

 \subsection{Certification of the results of the direct solver for higher frequencies} 
 
 Once the existence of the \emph{pre-asymptotic regime} is established, 
 and thus the efficiency of  the hierarchical representation of the matrix  to reduce computational times and memory requirements,
is demonstrated,  this subsection is devoted to the study of the accuracy of the method in this regime.

In Table~\ref{error_dens_single}, we report the values of the estimator  
  $\mathcal{I}(\delta_{\mathcal{H},F},\delta)$ and the observed normalized residual  $\frac{||{\bf b}-\mathbb{A}{\bf x}_0||_2}{ ||{\bf b} ||_2}$ if equation~(\ref{EFIE}) is solved with the $\mathcal{H}$-LU based direct solver.
 Four meshes are considered  and the density of points per S-wavelength is fixed consistently with the previous subsection.
 When no compression is possible ($\tau(\mathcal{H})=1$ for $N=1 \ 926$) the estimated and observed residuals are, as expected, of the order of the machine accuracy. 
 When the required accuracy is reduced by two orders of magnitude, the two residuals follow the same decrease. 
 The estimator computed in Frobenius norm is again not very sharp and the solver  could appear to be less accurate when the frequency increases.  But the discrepancy  between the true system matrix and its $\mathcal{H}$-matrix representation is seen to be of the order of $\varepsilon_{\operatorname{ACA}}$ even though computed in Frobenius norm. The increase of the estimator is due to the increase of the ratio $\frac{||{\bf x}_0||_2}{||{\bf b}||_2}$. In addition, the ratio   $\frac{\delta}{||{\bf b}||_2}$ indicates the very good stability of the $\mathcal{H}$-LU factorization.

\begin{table}[!htb]
  \begin{center}
 \begin{tabular}{ccc|cc|ccc}
 \hline
 $\#$ DOFs & $\omega$ & $\varepsilon_{\operatorname{ACA}}=\varepsilon_{\operatorname{LU}}$ & $\mathcal{I}(\delta_{\mathcal{H},F},\delta)$ & $\frac{||{\bf b}-\mathbb{A}{\bf x}_0||_2}{ ||{\bf b} ||_2}$ & $\delta_{\mathcal{H},F}$ & $\frac{||{\bf x}_0||_2}{||{\bf b}||_2}$ & $\frac{\delta}{||{\bf b}||_2}$ \\
 \hline
 $1 \ 926$ & $4$ & $10^{-4}$ & $3.12 \ 10^{-15}$ & $3.42 \ 10^{-15}$& $9.07 \ 10^{-17}$ & $8.37$ & $2.36 \ 10^{-15}$\\
 $7 \ 686$ & $8.25$ & $10^{-4}$ & $1.73 \ 10^{-3}$ & $4.28 \ 10^{-5}$ & $1.02 \ 10^{-4}$ & $16.83$ & $2.23 \ 10^{-5}$\\
 $30 \ 726$ & $16.5$ & $10^{-4}$ &$6.42 \ 10^{-3}$ & $1.68 \ 10^{-4}$ & $1.94 \ 10^{-4}$ & $32.66$ & $7.04 \ 10^{-5}$ \\
 $1 \ 926$ & $4$ & $10^{-6}$ & $3.12 \ 10^{-15}$ & $3.42 \ 10^{-15}$& $9.07 \ 10^{-17}$ & $8.37$ & $2.36 \ 10^{-15}$\\
 $7 \ 686$ & $8.25$ & $10^{-6}$ &$1.79 \ 10^{-5}$ & $2.61 \ 10^{-7}$ & $1.06 \ 10^{-6}$ & $16.83$ & $9.63 \ 10^{-8}$\\  
  $30 \ 726$ & $16.5$ & $10^{-6}$ & $7.14 \ 10^{-5}$ &$1.36 \ 10^{-6} $ & $2.18 \ 10^{-6}$& $32.66$&  $2.93 \ 10^{-7}$\\   
  \hline
 \end{tabular}
 \caption{Comparison of the estimation of the error   and the actual error of the direct solver for a fixed density of ten points per S-wavelength  for the single-layer potential.}
 \label{error_dens_single}
 \end{center}
 \end{table}

These results confirm the numerical efficiency and accuracy of $\mathcal{H}$-matrix based solvers even for higher frequencies. The estimator $\mathcal{I}(\delta_{\mathcal{H},F},\delta)$ gives a rough idea of the accuracy of the solver. But it is better to check separately the two components of the error.

 {
\begin{remark}
The computational times of the iterative and direct solvers have not been compared due to the difficulty to perform a fair and exhaustive comparison. On the one hand, the iterative solver can easily be optimized and parallelized but its performances depend largely on the quality of the preconditioner used to reduce the number of iterations. On the other hand, the optimization of the direct solver is a more involved operation that requires a task-based parallelization. Our experience with the two solvers is that the computational time of the direct solver is increasing very rapidly as a function of the number of DOFs if no parallelization is performed. However, the preconditioning of iterative BEM solvers for oscillatory problems is still an open question. 
\end{remark}
}

\section{Conclusion and future work}
 
 In this work, $\mathcal{H}$-matrix based iterative and direct solvers have been successfully extended to  3D elastodynamics in the frequency domain. 
  To perform the low-rank approximations in an efficient way, the Adaptive Cross Approximation has been modified to consider problems with vector unknowns.
  In addition,   the behavior of the hierarchical algorithm in the context of  oscillatory kernels has been carefully studied. It has been shown theoretically and then numerically confirmed  that the maximum numerical rank among all the $\eta$-\emph{admissible} blocks is constant if the frequency is fixed and the density of discretization points per wavelength is increased (i.e. for the low-frequency regime).
On the other hand, if the density of points per S-wavelength is fixed, as it is usually the case for higher frequency  problems,  the growth of the maximum numerical rank among all the $\eta$-\emph{admissible} blocks is limited to $O(N^{1/2})$.
Then, since the $\mathcal{H}$-LU factorization is based on various levels  of approximations, 
   a simple estimator to certify the results of the $\mathcal{H}$-LU direct solver has been derived.  It has been shown that not only the accuracy of the    $\mathcal{H}$-matrix representation can be monitored by a simple parameter but also the additional error introduced by the  $\mathcal{H}$-LU direct solver. 
Globally,  when  combined with the Boundary Element Method formulation, the $\mathcal{H}$-matrix representation of the system matrix permits to drastically reduce computational costs in terms of memory requirements and computational times for the frequency range that can be  considered with the current available computational resources.

Applications of the present $\mathcal{H}$-matrix accelerated BEM to 3D realistic cases in seismology are underway. Moreover, a natural extension of this work will be to
determine the efficiency of this approach for  isotropic and anisotropic visco-elastic  problems. 

\section*{Acknowledgments }
The authors would like to thank Shell for   funding the PhD of Luca Desiderio.

\section*{References}

\small
\bibliographystyle{plain}
\bibliography{Hmat}

\end{document}